\documentclass{amsart}
\usepackage[T1]{fontenc}
\usepackage{bussproofs}
\usepackage{amssymb}
\usepackage{mathabx}
\usepackage{hyperref}

\newtheorem{theorem}[subsection]{Theorem}
\newtheorem{lemma}[subsection]{Lemma}
\newtheorem{corollary}[subsection]{Corollary}

\newtheorem{remark}{Remark}
\newtheorem{example}{Example}
\newtheorem{definition}[subsection]{Definition}

\title[Approximate Completeness for GŁ$\forall$]{Approximate Completeness of Hypersequent Calculus for First-Order Łukasiewicz Logic}
\author{Jin Wei}
\date{}

\begin{document}
\begin{abstract}
Hypersequent calculus GŁ$\forall$ for first-order Łukasiewicz logic was first introduced by Baaz and Metcalfe, along with a proof of its approximate completeness with respect to standard $[0,1]$-semantics. The completeness result was later pointed out by Gerasimov that it only applies to prenex formulas. In this paper, we will present our proof of approximate completeness of GŁ$\forall$ for arbitrary first-order formulas by generalizing the original completeness proof to hypersequents. 
\end{abstract}

\maketitle

\section{Introduction}
Łukasiewicz logic is one of the fuzzy logics where formulas take value in $[0,1]$. It can be axiomatized in a Hilbert-style proof system that is sound and complete with respect to $[0,1]$-semantics, or arbitrary MV-algebras. Its first-order variant however takes the following approximate form:
$$A \text{ is valid in }[0,1]\iff A\oplus A^n \text{ is provable for all natural number }n$$
where $A\oplus A^n$ is a formula that approximates any truth assignment of $A$. The result is optimal in some sense as first-order valid formulas over $[0,1]$ are not recursively enumerable (see \cite{hajek}).

Gentzen-style proof systems for propostional Łukasiewicz logic were investigated by Metcalfe, Olivetti, and Gabbay in \cite{prop} in which a system of analytic hypersequent calculus GŁ was introduced. They prove the soundness and completeness of GŁ. A syntactical proof of cut elimination was also provided in \cite{synt}.

Baaz and Metcalfe extended the result to first-order Łukasiewicz logic in \cite{main} by introducing GŁ$\forall$ and proving its soundness and approximate completeness. The proof utilizes Skolemization and an approximate form of Herbrand's theorem to convert any given valid first-order formula $A$ to a valid quantifier-free formula $B\oplus B^n$. Then by the propositional completeness of GŁ, there exists a deduction of $B\oplus B^n$, from which a deduction of $A\oplus A^n$ can be found by re-introducing quantifiers. 

Later, Gerasimov points out  in \cite{gera} that this proof can only apply to prenex formulas, since we cannot prove equivalence between first-order formulas to its prenex form from GŁ$\forall$ without $cut$. It thus remains open whether GŁ$\forall$ is complete with respect to $[0,1]$-semantics for arbitrary first-order formulas or not.

In this paper, we provide an affirmative answer to the question by providing a proof for the approximate completeness of GŁ$\forall$. In fact, we have a stronger result:
$$\mathcal{H} \text{ is a valid hypersequent in }[0,1]\iff \mathcal{H}_{\frac{1}{n}} \text{ is provable in GŁ$\forall$ for all }n$$
where $\mathcal{H}_{\frac{1}{n}}$ is a hypersequent that approximates any truth interpretation of $\mathcal{H}$. The original statement in \cite{main} now becomes a corollary of this result.

\section{Preliminary}
\subsection{Propositional Łukasiewicz Logic}
The language of Łukasiewicz logic contains two connectives $\to$ and $\bot$. We have a Hilbert-style proof system HŁ for Łukasiewicz logic.
\begin{definition}[HŁ]
HŁ consists of four axioms and \textit{modus ponens}. 
 
$\vdash A\to (B\to A)$

$\vdash (A\to B)\to ((B\to C)\to (A\to C))$

$\vdash ((A\to B)\to B)\to ((B\to A)\to A)$ 

$\vdash (\neg B\to \neg A)\to (A\to B)$

\begin{prooftree}
    \AxiomC{$\vdash A$}
    \AxiomC{$\vdash A\to B$}
    \RightLabel{$mp$}
    \BinaryInfC{$\vdash B$}
\end{prooftree}
\end{definition}

\begin{definition}
A interpretation $v$ is a function on the set of formulas such that: 
\begin{itemize}
    \item $v(A) \in [0,1]$ for atomic formula $A$.
    \item $v(\bot):=1$.
    \item $v(A\to B):=v(B)\dotdiv  v(A)$.
\end{itemize}

We say a propositional formula $A$ is valid if $v(A)=0$ for all interpretation $v$ and we denote it as $\models_{[0,1]}A$.
\end{definition}

\begin{definition}
We provide a number of connectives definable using $\bot$ and $\to$.
\begin{itemize}
	\item $\neg A := A \to \bot$ and $\top :=\neg\bot$ 
	\item $A \oplus B:=\neg A \to B$ and $A \odot B :=\neg (\neg A \oplus \neg B)$
	\item $A \vee B:=(A\to B)\to B$ and $A\wedge B:= \neg (\neg A \vee \neg B)$
	\item $0.A:=\bot$ and $(n+1).A:=A\oplus n.A$
	\item $A^0:=\top$ and $A^{n+1}:=A\odot A^n$ 
\end{itemize}

\end{definition}

The interpretations and properties of the above connectives can be found in \cite{hajek} and \cite{main}.

\begin{remark}
Note that some authors use $0$ as the value of $\bot$ and $1$ as the value. The difference should be purely notational and our choice has certain advantages as it appeals to continuous model theorists and makes many notations simpler in this paper. 
\end{remark}

Towards a Gentzen-style proof system of Łukasiewicz logic, we first define a notation of sequents and hypersequents.

\begin{definition}
A sequent is of the form $\Gamma\Rightarrow \Delta$ where $\Gamma$ and $\Delta$ are finite multisets of formulas. A hypersequent is a finite multisets of sequents, denoted as $\Gamma_1\Rightarrow \Delta_1\mid \cdots \mid \Gamma_n\Rightarrow\Delta_n$.
\end{definition}

Note that different from classical logics, sequents and hypersequents in Łukasiewicz logic are based on multisets of formulas due to its connection to affine linear logic see (\cite{prop}). 

\begin{definition}
We extend valuation $v$ to sequents and hypersequents, whose values are in $\mathbb{R}$. Notice that here we use regular sum $\sum$ and subtraction $-$ instead of bounded ones.

\begin{itemize}
    \item $v(\Gamma\Rightarrow \Delta):=\sum_{A\in 
    \Delta}v(A)- \sum_{B\in 
    \Gamma}v(B)$.
    \item $v(\Gamma_1\Rightarrow \Delta_1\mid \cdots \mid \Gamma_n\Rightarrow\Delta_n):=\min_{1\leq i\leq n} v(\Gamma_i\Rightarrow \Delta_i)$

\end{itemize}

We say a hypersequent $\mathcal{H}$ with no free variables is valid if $v(\mathcal{H})\leq 0$ for all valuation $v$ and we denote as $\models_{[0,1]}\mathcal{H}$. 

\end{definition}

Informally, $\mathcal{H}$ is valid if and only if some of its component represents an inequality valid in all interpretations. Now we are ready to define the hypersequent calculus GŁ for propositional Łukasiewicz logic.

\begin{definition}[GŁ]
A hypersequent $\mathcal{H}$ is provable from GŁ, denoted as $\vdash_{\text{GŁ}} \mathcal{H}$, if it is derivable from the following rules. 
\vspace{0.3cm}

Initial sequents:

\begin{center}
    \AxiomC{}
    \RightLabel{id $\qquad\qquad\qquad\quad$}
  \UnaryInfC{$A \Rightarrow A$}
 \DisplayProof 
     \AxiomC{}
    \RightLabel{$\Rightarrow\qquad\qquad\qquad\quad$}
  \UnaryInfC{$\Rightarrow $}
 \DisplayProof 
    \AxiomC{}
    \RightLabel{$\bot\Rightarrow$}
  \UnaryInfC{$\bot\Rightarrow A$}
 \DisplayProof 
 \end{center}
\vspace{0.3cm}

Structural rules:

\begin{center}
%     \AxiomC{$\mathcal{G}\mid \mathcal{H}$}
%     \RightLabel{ex$\qquad\qquad\qquad\qquad\quad\qquad$}
% \UnaryInfC{$\mathcal{H}\mid \mathcal{G}$}
%  \DisplayProof 
    \AxiomC{$\mathcal{G}\mid \Gamma\Rightarrow \Delta\mid \Gamma\Rightarrow \Delta$}
    \RightLabel{ec$\qquad\qquad\qquad\qquad$}
\UnaryInfC{$\mathcal{G}\mid \Gamma\Rightarrow \Delta$}
 \DisplayProof 
     \AxiomC{$\mathcal{G}\mid \Gamma\Rightarrow \Delta$}
    \RightLabel{ew}
    \UnaryInfC{$\mathcal{G}\mid \Gamma\Rightarrow \Delta
    \mid \Gamma'\Rightarrow \Delta'$}
   \DisplayProof 
\end{center}

\vspace{0.3cm}

\begin{center}
     \AxiomC{$\mathcal{G}\mid \Gamma_0,\Gamma_1\Rightarrow \Delta_0,\Delta_1$}
    \RightLabel{split$\qquad\qquad$}
    \UnaryInfC{$\mathcal{G}\mid \Gamma_0\Rightarrow \Delta_0\mid \Gamma_1\Rightarrow \Delta_1$}
\DisplayProof
    \AxiomC{$\mathcal{G}\mid \Gamma_0\Rightarrow \Delta_0$}
    \AxiomC{$\mathcal{G}\mid \Gamma_1\Rightarrow \Delta_1$}
    \RightLabel{mix}
\BinaryInfC{$\mathcal{G}\mid \Gamma_0,\Gamma_1\Rightarrow \Delta_0,\Delta_1$}
\DisplayProof
\end{center}

\vspace{0.3cm}

\begin{center}

     \AxiomC{$\mathcal{G}\mid \Gamma\Rightarrow \Delta$}
    \RightLabel{wl}
    \UnaryInfC{$\mathcal{G}\mid \Gamma,A\Rightarrow \Delta$}
 \DisplayProof
\end{center}

\vspace{0.3cm}

Logical rules: 

\begin{center}
    \AxiomC{$\mathcal{G}\mid \Gamma,B\Rightarrow A, \Delta \mid \Gamma\Rightarrow \Delta$}
    \RightLabel{$\to \Rightarrow\qquad$}    
    \UnaryInfC{$\mathcal{G}\mid \Gamma,A\to B\Rightarrow \Delta$}
\DisplayProof
      \AxiomC{$\mathcal{G}\mid \Gamma,A\Rightarrow B,\Delta$}
   \AxiomC{$\mathcal{G}\mid \Gamma\Rightarrow \Delta$}
    \RightLabel{$\Rightarrow \to $}
    \BinaryInfC{$\mathcal{G}\mid \Gamma\Rightarrow A\to B,\Delta$}
\DisplayProof
\end{center}
\end{definition}

\begin{definition}[Cut]
We have the following form of $cut$ rule in addition to GŁ.
\begin{prooftree}
     \AxiomC{$\mathcal{G}\mid \Gamma,A\Rightarrow A,\Delta$}
     \RightLabel{$cut$}
    \UnaryInfC{$\mathcal{G}\mid \Gamma\Rightarrow \Delta$}
\end{prooftree}
\end{definition}

\begin{theorem}
For any formula $A$, the following are equivalent
\begin{enumerate}
    \item $\vdash_{\text{HŁ}}A$
    \item $\vdash_{\text{GŁ} + cut}\Rightarrow A$
    \item $\vdash_{\text{GŁ}}\Rightarrow A$
    \item $\models_{[0,1]} A$
\end{enumerate}
\begin{proof}
See proofs in \cite{prop} and \cite{synt}. 
\end{proof}
\end{theorem}

\subsection{First-Order Łukasiewicz Logic} We extend the propositional Łukasiewicz logic to first-order by adding existential quantifier $\exists$. 

\begin{definition}

GŁ$\forall$ is obtained from GŁ by adding the following two rules:
\begin{center}
	    \AxiomC{$\mathcal{G}\mid \Gamma,A(c)\vdash \Delta$}
    \RightLabel{$\exists\Rightarrow\qquad\qquad\qquad\qquad$}
    \UnaryInfC{$\mathcal{G}\mid \Gamma,(\exists x)A[x/c]\Rightarrow \Delta$}
\DisplayProof
      \AxiomC{$\mathcal{G}\mid \Gamma\Rightarrow A[t/x],\Delta$}
    \RightLabel{$\Rightarrow\exists$}
    \UnaryInfC{$\mathcal{G}\mid \Gamma\Rightarrow (\exists x)A,\Delta$}
\DisplayProof
\end{center}
where $c$ does not appear free elsewhere in the hypersequent and $t$ is an arbitrary term.
\end{definition}

Universal quantification $(\forall x)A$ is understood as a shorthand for $\neg ((\exists x)\neg A)$ (slightly inconsistent with the name GŁ$\forall$). Now we introduce the semantics for first-order Łukasiewicz logic. 

\begin{definition}
A $[0,1]$-structure $\mathfrak{M}$ is a set $M$ together with its interpretation for function and relation symbols. 
\begin{itemize}
	\item If $f$ is an $n$-ary function symbol, then $f^\mathfrak{M}$ is a function $M^n\to M$. 
	\item If $R$ is an $n$-ary relation symbol, then $R^\mathfrak{M}$ is a function $M^n\to [0,1]$.
   \item $\bot^{\mathfrak{M}}=1$.
   \item $(A\to B)^{\mathfrak{M}}=B^{\mathfrak{M}}\dotdiv A^{\mathfrak{M}}$.
   \item $((\exists x) A)^{\mathfrak{M}}=\inf_{m\in M} (A[m/x])^{\mathfrak{M}}$.
    \item $(\Gamma\Rightarrow \Delta)^{\mathfrak{M}}:=\sum_{A\in 
    \Delta}A^{\mathfrak{M}}- \sum_{B\in 
    \Gamma}B^{\mathfrak{M}}$.
    \item $(\Gamma_1\Rightarrow \Delta_1\mid \cdots \mid \Gamma_n\Rightarrow\Delta_n)^{\mathfrak{M}}:=\min_{1\leq i\leq n} (\Gamma_i\Rightarrow \Delta_i)^{\mathfrak{M}}$
\end{itemize}

We say $\mathcal{H}$ is valid if $\mathcal{H}$ contains no free variable and $\mathcal{H}^{\mathfrak{M}}\leq 0$ for all $[0,1]$-structures $\mathfrak{M}$. We denote it as $\models_{[0,1]}\mathcal{H}$.
\end{definition}

\begin{definition}
We adopt the following notations for any given formula $A$, multisets $\Gamma,\Delta$, and natural number $n$:
\begin{itemize}
    \item $nA$ denotes the multiset $\{\overbrace{A,\ldots,A}^n\}$
    \item $n\Gamma$ denotes the multiset $\overbrace{\Gamma\cup \cdots \cup \Gamma}^n$.
    \item $\Gamma\Rightarrow_{ \frac{1}{n}}\Delta$ denotes the sequent $$\bot, n\Gamma\Rightarrow n\Delta$$
    whose intended meaning is that $\Gamma\Rightarrow \Delta$ is valid up to an $\frac{1}{n}$ error in all interpretations, since $\bot$ has truth value $1$.
    \item We generalize the above notation to hypersequents $\mathcal{H}$. Say $\mathcal{H}$ is of the form $\Gamma_1\Rightarrow\Delta_1\mid \cdots \mid \Gamma_k\Rightarrow \Delta_k$. Then $\mathcal{H}_{\frac{1}{n}}$ is the hypersequent $$\Gamma_1\Rightarrow_{\frac{1}{n}}\Delta_1\mid \cdots \mid \Gamma_k\Rightarrow_{\frac{1}{n}} \Delta_k$$
    \item $A[\bar t/\bar x], \Gamma[\bar t/\bar x], \mathcal{H}[\bar t/\bar x]$ represent substitutions of free variables $\bar x$ by terms $\bar t$ on formula, sequent, and hypersequent levels.

    \item A formula is atomic if it is either a propositional variable/predicate or $\bot$.
\end{itemize}
\end{definition}

\section{Completeness of GŁ}
We independently developed a proof of completeness of GŁ with respect to $[0,1]$-semantics, unaware that similar ideas are already hinted in \cite{prop} and \cite{alte}. Readers may skip to Section \ref{GŁforall} for the main result of the paper. 

\begin{lemma}\label{atom}
If $\Gamma\Rightarrow \Delta$ is a valid sequent consisting only of atomic formulas, then we have a deduction of $\Gamma\Rightarrow \Delta$ from GŁ.
\begin{proof}
If the sequent is of the form $\Gamma,A\Rightarrow A,\Delta$, then $\Gamma\Rightarrow \Delta$ remains valid. Furthermore, it suffices to derive $\Gamma\Rightarrow \Delta$ by the following:
\begin{prooftree}
    \AxiomC{$\Gamma\Rightarrow \Delta$}
    \AxiomC{}
    \RightLabel{$id$}
    \UnaryInfC{$A \Rightarrow A$}
    \RightLabel{$mix$}
    \BinaryInfC{$\Gamma,A \Rightarrow A,\Delta$}
\end{prooftree}

Now we may assume without loss of generality that $\Gamma\cap \Delta=\emptyset$. 

If $\Delta$ is empty, we are done either by $wl$ if $\Gamma$ is not empty or $\Rightarrow$ if $\Gamma$ is. Otherwise, the sequent is now the form $\Gamma\Rightarrow n_1A_1, \ldots,n_k A_k$. Suppose there are $n<n_1+\cdots + n_k$-many $\bot$ in $\Gamma$, we have an interpretation $v$ where $v(A_i)=1$ for each $i$ and $v(B)=0$ otherwise. Then $v(\Gamma\Rightarrow n_1A_1, \ldots,n_k A_k)<n_1+\cdots + n_k-n>0$, which contradicts with validity of $\Gamma \Rightarrow n_1A_1, \ldots,n_k A_k$. Therefore, $n\geq n_1+\cdots + n_k$ and we now construct a deduction of $\Gamma\Rightarrow\Delta$ by $n_1+\cdots + n_k$-many $\bot\Rightarrow$ rules together with $wl$ and $mix$.
\end{proof}
\end{lemma}

We will prove completeness of GŁ by a bottom-up proof search. We need the following variant of Farkas' lemma from linear programming to construct the leafs hypersequents. The proof of the variant can be found in \cite{alt2}. 

\begin{lemma}\label{thmofalt}
Given matrices $M\in \mathbb{R}^{m_1\times n}$ and $N\in \mathbb{R}^{m_2\times n}$ and vectors $a\in \mathbb{R}^{m_1}$ and $b\in \mathbb{R}^{m_2}$. Exactly one of the following holds.
\begin{itemize}
    \item There exists $x\in \mathbb{R}^n$ satisfying 
    \begin{eqnarray*}
        M x&\leq& a \\
        N x&\ll& b\\
        x &\geq& 0
    \end{eqnarray*}
    \item There exist $\lambda\in \mathbb{R}^{m_1}$ and $\mu\in \mathbb{R}^{m_2}$ satisfying 
    \begin{eqnarray*}
        \lambda &\geq & 0 \\
        \mu &\geq & 0 \\
        \lambda^T M + \mu^T N &\geq& 0\\
        \lambda^T a + \mu^T b &\leq& 0\\
        \lambda^T a + \mu^T b < 0 \text{ or } \mu &>&0
    \end{eqnarray*}
\end{itemize}
where $v\leq w$ means $v_i\leq w_i$ for all $i$ and $v\ll w$ means $v_i<w_i$ for all $i$.
\end{lemma}

\begin{theorem}[Completeness of GŁ]\label{GŁcmp}
$$\text{If } \mathcal{H} \text{ is a valid hypersequent, then } \mathcal{H}\text{ is derivable from GŁ.}$$
\begin{proof}
We construct a search tree for a deduction of valid hypersequent $\mathcal{H}$.

\begin{itemize}
    \item $\mathcal{H}$ is the base node.
    \item At node $\mathcal{G}\mid \Gamma \Rightarrow \Delta, A\to B$, search for deductions of $\mathcal{G}\mid \Gamma ,A \Rightarrow \Delta, B$ and $\mathcal{G}\mid \Gamma \Rightarrow \Delta$.
    \item At node $\mathcal{G}\mid \Gamma, A\to B \Rightarrow \Delta$, search for deduction of $\mathcal{G}\mid \Gamma ,A \Rightarrow \Delta, B\mid \Gamma\Rightarrow \Delta$. 
    \item At node $\mathcal{G}\mid \Gamma,A \Rightarrow A,\Delta$ where $A$ is either $\bot$ or atomic, search for a deduction of $\mathcal{G}\mid \Gamma \Rightarrow \Delta$. 
\end{itemize}

Therefore, it suffices to deal with hypersequents consisting of only atomic formulas. Let $A_1,\ldots, A_m$ enumerate all propositional variables appearing in $\mathcal{H}$ and say $\mathcal{H}=\Gamma_1\Rightarrow \Delta_1\mid \ldots\mid \Gamma_k\Rightarrow \Delta_k$. Since $\mathcal{H}$ is a valid hypersequent, we know that the following system of inequalities of $m\times m$-matrix $M$ and $k\times m$-matrix $N$
\begin{eqnarray*}
  v(A_1) &\leq& 1 \\
 &\vdots& \\
 v(A_m) &\leq& 1   \\
   -a_{11} v(A_1) -\ldots - a_{1m} v(A_m) &<& -c_0   \\
   \vdots \\
  -a_{k1} v(A_1) -\ldots - a_{km} v(A_m) &<& -c_k
\end{eqnarray*}
has no nonnegative solution for any truth-interpretation $v$ where $a_{ij}$ represents the multiplicity of $A_j$ in $\Gamma_i\Rightarrow \Delta_i$ where $a_{ij}$ is negative if $x_j\in \Gamma_i$, positive if $x_j\in \Delta_i$, and zero otherwise; $c_i$ represents the multiplicity of $\bot$ in $\Gamma_i\Rightarrow \Delta_i$ with the same polarity. 

Applying Lemma \ref{thmofalt} to the above system of inequalities, we obtain $\lambda\geq 0$ and $\mu\geq 0$. We may further assume $\lambda$ and $\mu$ are vectors of non-negative integers since all coefficients are integers. We obtain the following inequalities:
\begin{eqnarray*}
    \lambda_j -\sum_{i=1}^k \mu_i a_{ij} \geq 0 \text{ for each $i$}\\ 
    -\sum_{j=1}^m \lambda_j +\sum_{i=1}^k \mu_i c_i \geq 0
\end{eqnarray*}

Hence the following inequality holds for all $v$:
$$ (\sum_{j=1}^m \lambda_j-\sum_{i=1}^k \mu_i c_i)\cdot 1+(-\lambda_1 +\sum_{i=1}^k \mu_i a_{i1} )\cdot v(A_1) +\cdots +( -\lambda_m +\sum_{i=1}^k \mu_i a_{im} )\cdot v(A_m)\leq 0$$

Combining with $\lambda_i$-many inequalities $v(A_i)\leq 1$ for each $i$, we simplify it into:
$$ (\sum_{i=1}^k \mu_i c_i)\cdot 1+(\sum_{i=1}^k \mu_i a_{i1} )\cdot v(A_1) +\cdots +(\sum_{i=1}^k \mu_i a_{im} )\cdot v(A_m)\leq 0$$

Construct a sequent $\Gamma \Rightarrow \Delta$ where we add $\mu_ia_{ij}$-many $A_j$ in $\Gamma$ if $a_{ij}$ is positive $t$ and to $\Delta$ if $a_{ij}$ is negative; similarly, we add $\bot$ according to polarity of $c_i$. By design, $\Gamma \Rightarrow \Delta$ is a valid sequent consisting only of atomic formulas since the above inequality holds for all interpretation $v$. We then have a deduction of $\Gamma\Rightarrow \Delta$ by Lemma \ref{atom}. We can recover the original $\mathcal{H}$ by the following derivation where the last step is a $ew$ if $\mu_i=0$ and a $ec$ if $\mu_i\neq 0$. 

\begin{prooftree}
    \AxiomC{$\Gamma\Rightarrow \Delta$}
    \RightLabel{$split$}
    \UnaryInfC{$\underbrace{\Gamma_1\Rightarrow \Delta_1 \mid \cdots \mid \Gamma_1\Rightarrow \Delta_1}_{\mu_1}\mid \cdots \mid \underbrace{\Gamma_k\Rightarrow \Delta_k \mid \cdots \mid \Gamma_k\Rightarrow \Delta_k}_{\mu_k}$}
    \RightLabel{$ec/ew$}
    \UnaryInfC{$\Gamma_1\Rightarrow\Delta_1\mid \cdots \mid \Gamma_k\Rightarrow\Delta_k$}
\end{prooftree}
\end{proof}
\end{theorem}

\section{Completeness of GŁ$\forall$}\label{GŁforall}
Given a valid hypersequent $\mathcal{H}$ in first-order Łukasiewicz logic (in particular, $\mathcal{H}$ contains no free variables), we will construct a search tree for deductions of $\mathcal{H}$. The following formula serves as a motivating example for the actual search algorithm. 

\begin{example}[Drinker Paradox] Define 
$\varphi:=\exists x (A(x)\to \forall y A(y))$ where $A$ is a unary predicate symbol.

With some computation, we have that $\models_{[0,1]}\varphi$ and $\vdash_{\text{GŁ$\forall$}+cut} \Rightarrow\varphi$. However, $\varphi$ is not provable in GŁ$\forall$ (see \cite{main}). Instead, we do have derivations of $\Rightarrow_{\frac{1}{n}}\varphi$ for each $n$. To search for such derivations, we do a bottom-up  Skolemization. 

\begin{prooftree}
    \AxiomC{$A(x)\Rightarrow A(f(x))$}
    \UnaryInfC{$A(x),\bot \Rightarrow  A(f(y)),\bot $}
    \UnaryInfC{$A(x),A(f(y))\to \bot\Rightarrow  \bot $}
    \UnaryInfC{$A(x),\exists y (A(y)\to \bot )\Rightarrow  \bot $}
    \UnaryInfC{$A(x)\Rightarrow  \exists y (A(y)\to \bot )\to \bot $}
    \UnaryInfC{$ A(x)\Rightarrow\forall y A(y)$}
    \UnaryInfC{$\Rightarrow A(x)\to \forall y A(y)$}
    \UnaryInfC{$\Rightarrow\exists x (A(x)\to \forall y A(y))$}
\end{prooftree}

The sequent $A(x)\Rightarrow A(f(x))$ may not have any Herbrand witnesses. Instead, for each $n$, there are approximate Herbrand terms $x,f(x),f(f(x))\ldots,f^{n-1}(x)$ such that the following is valid
$$A(x)\Rightarrow_{\frac{1}{n}} A(f(x))\mid A(f(x))\Rightarrow_{\frac{1}{n}} A(f^2(x))\mid \cdots \mid A(f^{n-1}(x))\Rightarrow_{\frac{1}{n}} A(f^n(x))$$
We guarantee to have a derivation of the above quantifier-free hypersequent by Completeness of GŁ. 

\begin{prooftree}
\AxiomC{}
\UnaryInfC{$n\bot\Rightarrow nA(f^n(x))$}
\UnaryInfC{$n\bot, nA(x)\Rightarrow nA(f^n(x))$}
\UnaryInfC{$n\bot, nA(x),nA(f(x)),\ldots,nA(f^{n-1}(x))\Rightarrow nA(f(x)),\ldots,nA(f^n(x))$}
\UnaryInfC{$\bot, nA(x)\Rightarrow nA(f(x))\mid \bot, nA(f(x))\Rightarrow nA(f^2(x))\mid \cdots \mid \bot, nA(f^{n-1}(x))\Rightarrow n A(f^n(x))$}
\UnaryInfC{$A(x)\Rightarrow_{\frac{1}{n}} A(f(x))\mid A(f(x))\Rightarrow_{\frac{1}{n}} A(f^2(x))\mid \cdots \mid A(f^{n-1}(x))\Rightarrow_{\frac{1}{n}} A(f^n(x))$}
\end{prooftree}
Then we will follow the above search tree from top-down to recover $\Rightarrow_{\frac{1}{n}}\exists x(A(x)\to \forall y A(y)$ by re-introducing connectives in a certain order. We also need to admit $\frac{1}{n}$-version of connective introduction rules for that. 
\end{example}

The above example can be generalized to arbitrary valid hypersequents to yield a completeness of GŁ$\forall$ over $[0,1]$-semantics. 

\begin{definition}[Skolemization Tree]\label{dctr}
Given a hypersequent $\mathcal{H}$, we build a binary tree $T_{\mathcal{H}}$ of hypersequents with free variables $\bar x$ to remove quantifiers and $\to$. The order can be taken arbitrarily if multiple conditions of the following are satisfied.

Base node:
$\mathcal{H}$ is the base node.
    
\textbf{Right $\to$ case:}
If the left-most leaf $\mathcal{G}$ has the form $$\mathcal{L}\mid \Gamma\Rightarrow A\to B,\Delta$$ 
    we attach to $\mathcal{G}$ a left node $\mathcal{L}\mid \Gamma, A\Rightarrow B,\Delta$ and a right node $\mathcal{L}\mid \Gamma\Rightarrow \Delta$. We say the position of sequent $ \Gamma\Rightarrow \Delta$ on the right node corresponds to the position of sequent $ \Gamma,A\Rightarrow B,\Delta$ on the left node (we may take hypersequents as lists of sequents with $exchange$ rule if necessary). 
    
    Furthermore, for any other leaf $\mathcal{R}$, suppose $\mathcal{R}$ has the form $\mathcal{L}^*\mid \Gamma^*\Rightarrow A\to B, \Delta^*$ where the position of $\Gamma^*\Rightarrow A\to B, \Delta^*$ in $\mathcal{R}$ corresponds to $\Gamma\Rightarrow A\to B, \Delta$ in $\mathcal{G}$. We attach to $\mathcal{R}$ two nodes $\mathcal{L}^*\mid \Gamma^*, A\Rightarrow B,\Delta^*$ and $\mathcal{L}^*\mid \Gamma^*\Rightarrow \Delta^*$. Otherwise, $A\to B$ is absent in the position and we attach $\mathcal{R}$ to itself in that case. 

\textbf{Left $\to$ case:}
If the left-most leaf $\mathcal{G}$ has the form $$\mathcal{L}\mid \Gamma,A\to B\Rightarrow \Delta$$
    we attach to $\mathcal{G}$ the node $\mathcal{L}\mid \Gamma, B\Rightarrow A,\Delta\mid \Gamma\Rightarrow \Delta$. We assume $\Gamma, B\Rightarrow A,\Delta$ take the position of $\Gamma, A\to B\Rightarrow \Delta$ and create a new position for $\Gamma\Rightarrow \Delta$.
    
    For any other leaf $\mathcal{R}$, suppose $\mathcal{R}$ has the form $ \mathcal{L}^*\mid \Gamma^*,A\to B\Rightarrow \Delta^*$ where the position of $ \Gamma^*,A\to B\Rightarrow \Delta^*$ in $\mathcal{R}$ corresponding to $\Gamma,A\to B\Rightarrow \Delta$ in $\mathcal{G}$. We attach to $\mathcal{R}$ the node $ \mathcal{L}^*\mid \Gamma^*,B\Rightarrow A, \Delta^*\mid \Gamma^*\Rightarrow \Delta^*$, while letting $\Gamma^*,B\Rightarrow A,\Delta^*$ correspond to $\Gamma,B\Rightarrow A,\Delta$ and $\Gamma^*\Rightarrow \Delta^*$ correspond to $\Gamma\Rightarrow \Delta$. Otherwise, $A\to B$ is absent in the position and we attach $\mathcal{R}$ to itself in that case.

\textbf{Right $\exists$ case:}
If the left-most leaf $\mathcal{G}$ takes the form $$\mathcal{L}\mid \Gamma\Rightarrow (\exists y) A,\Delta$$ 
    we attach to $\mathcal{G}$ the node $\mathcal{L}\mid\Gamma\Rightarrow  A[x/y],\Delta$ where $x$ is a new variable symbol. 
    
    For any other leaf $\mathcal{R}$, suppose $\mathcal{R}$ has the form $ \mathcal{L}^*\mid \Gamma^*\Rightarrow (\exists y) A,\Delta^*$ where the position of $\Gamma^*\Rightarrow (\exists y) A,\Delta^*$ in $\mathcal{R}$ corresponds to $\Gamma\Rightarrow (\exists y) A,\Delta$ in $\mathcal{G}$. We attach to $\mathcal{R}$ the node $\mathcal{L}^*\mid \Gamma^*\Rightarrow  A[x/y],\Delta^*$ with the same variable $x$. Otherwise, $(\exists y) A$ is absent in the position and we attach $\mathcal{R}$ to itself in that case.

\textbf{Left $\exists$ case:}
If the left-most leaf $\mathcal{R}$ takes the form $$\mathcal{L}\mid \Gamma, (\exists y)A\Rightarrow\Delta$$ 
    we attach to $\mathcal{G}$ the node $\mathcal{G}\mid \Gamma, A[f(\bar x)/y]\Rightarrow\Delta$ where $f$ is a new function symbol and $\bar x$ enumerates all free variables appearing at the current stage of the tree. 
    
    For all other leaf $\mathcal{R}$, suppose $\mathcal{R}$ has the form $ \mathcal{G}^*\mid \Gamma^*,(\exists y) A\Rightarrow \Delta^*$ where the position of $\Gamma^*,(\exists y) A\Rightarrow \Delta^*$ in $\mathcal{R}$ corresponds to $\Gamma,(\exists y) A\Rightarrow \Delta$ in $\mathcal{G}$. We attach to $\mathcal{R}$ the node $\mathcal{G}^*\mid \Gamma^*,A[f(\bar x)/y]\Rightarrow  \Delta^*$ with the same $f$ and $\bar x$. Otherwise, $(\exists y) A$ is absent in the position and we attach $\mathcal{R}$ to itself in that case.
\end{definition}

Let $h$ denote the max length of all paths in $T_{\mathcal{H}}$. We define level sets collecting all nodes sharing the same depth. 
\begin{definition}[Level Sets]
For each $i\leq h$, define $$X_i:=\{\sigma(i): \sigma \text{ is a path in }T_{\mathcal{H}}\}$$
\end{definition}

\begin{lemma}\label{sync}
For each hypersequent $\mathcal{H}$, the left-most path in $T_{\mathcal{H}}$ has length $h$ and ends with a hypersequent with no quantifiers. For each $i\leq h$, say $\mathcal{G}$ and $\mathcal{R}$ are in $X_i$ where $\mathcal{G}$ is in the left-most path. Then $\mathcal{R}$ has the same depth as $\mathcal{G}$ and any formula appearing in $\mathcal{R}$ also appears in $\mathcal{G}$.
\begin{proof}
We prove by an induction on $i$. The tree construction terminates exactly when the left-most leaf consists no $\to$ and quantifiers. Tree construction on other branches always follow the construction on the left-most branch, except missing some formulas due to the Right $\to$ case. 
\end{proof}
\end{lemma}  

Given a valid hypersequent $\mathcal{H}$, the level set $X_i$ of $T_{\mathcal{H}}$ keeps an approximate form of satisfiability. 

\begin{lemma}\label{apsa}
Say $\mathcal{H}$ is a valid hypersequent and $i\leq h$. Then for any structure $\mathfrak{M}$ and any natural number $n$, the following inequality holds
$$\inf_{\bar m\in M} \max_{\mathcal{G}\in X_i}(\mathcal{G}_{\frac{1}{n}}[\bar m/\bar x])^{\mathfrak{M}}\leq 0$$
where $\bar x$ are the free variables appearing in $X_i$.
\begin{proof}
We proceed by induction on $i$. 

%Base case

In the base case for $i=1$, we have that $X_1=\{\mathcal{H}\}$. Since $\mathcal{H}$ is valid and contains no free variable, we have for each $\mathfrak{M}$ and $n$ that:
$$\inf_{\bar m\in M} (\mathcal{H}_{\frac{1}{n}}[\bar m/\bar x])^{\mathfrak{M}}=(\mathcal{H}_{\frac{1}{n}})^{\mathfrak{M}}=n(\mathcal{H}^{\mathfrak{M}}-\frac{1}{n})< n\cdot \mathcal{H}^{\mathfrak{M}}\leq 0$$

In the inductive case for $i+1$, there are four cases based on the left-most node in $X_{i+1}$. Let $\mathfrak{M}$ and $n$ be arbitary.

%Right \to

\textbf{Right $\to$ case:}
When the left-most node in $X_{i+1}$ takes the form $\mathcal{L}\mid \Gamma,A\Rightarrow B,\Delta$ with parent node $\mathcal{L}\mid \Gamma \Rightarrow A\to B,\Delta$ in $X_i$, the induction hypothesis give
$$\inf_{\bar m}\max_{G\in X_i}(\mathcal{G}_{\frac{1}{n}}[\bar m/\bar x])^{\mathfrak{M}}\leq 0$$
It implies for each $\epsilon>0$ there are $\bar m\in M$ such that 
$$\max_{G\in X_i}(\mathcal{G}_{\frac{1}{n}}[\bar m/\bar x])^{\mathfrak{M}}< \epsilon$$
In particular, we have that
$$(\mathcal{L}_{\frac{1}{n}}\mid \Gamma\Rightarrow_{\frac{1}{n}} A\to B,\Delta)[\bar m/\bar x]^{\mathfrak{M}}< \epsilon$$
Since $(A\to B)^{\mathfrak{M}}=B^{\mathfrak{M}}\dotdiv A^{\mathfrak{M}}=\max(B^{\mathfrak{M}}- A^{\mathfrak{M}},0)$, we have 
$$    \max((\mathcal{L}_{\frac{1}{n}}\mid \Gamma,A\Rightarrow_{\frac{1}{n}}B,\Delta)[\bar m/\bar x]^{\mathfrak{M}},(\mathcal{L}_{\frac{1}{n}}\mid \Gamma\Rightarrow_{\frac{1}{n}}\Delta)[\bar m/\bar x]^{\mathfrak{M}})< \epsilon$$
and thus obtain two inequalities
\begin{eqnarray*}
    &&(\mathcal{L}_{\frac{1}{n}}\mid \Gamma,A\Rightarrow_{\frac{1}{n}}B,\Delta)[\bar m/\bar x]^{\mathfrak{M}}< \epsilon\\
    &&(\mathcal{L}_{\frac{1}{n}}\mid \Gamma\Rightarrow_{\frac{1}{n}}\Delta)[\bar m/\bar x]^{\mathfrak{M}}< \epsilon
\end{eqnarray*}
The first inequality covers $\mathcal{L}\mid \Gamma,A\Rightarrow B,\Delta$ in $X_{i+1}$ and the second inequality covers the element $\mathcal{L}\mid \Gamma\Rightarrow \Delta$ in $X_{i+1}$.

For the rest of elements in $X_{i+1}$, notice that there are two possibilities due to the definition of $T_{\mathcal{H}}$.
\begin{enumerate}
    \item Nodes $\mathcal{L}^*\mid \Gamma^*,A\Rightarrow B,\Delta^*$ and $\mathcal{L}^*\mid \Gamma^*\Rightarrow \Delta^*$ in $X_{i+1}$ have their parent node $\mathcal{L}^*\mid \Gamma^*\Rightarrow A\to B,\Delta^*$ in $X_i$. We then also have
    $$(\mathcal{L}^*_{\frac{1}{n}}\mid \Gamma^*\Rightarrow_{\frac{1}{n}}A\to B,\Delta^*)[\bar m/\bar x]^{\mathfrak{M}}<\max_{G\in X_i}(\mathcal{G}_{\frac{1}{n}}[\bar m/\bar x])^{\mathfrak{M}}<\epsilon$$
    In a similar proof, we derive the two desired inequalities $(\mathcal{L}^*_{\frac{1}{n}}\mid \Gamma^*,A\Rightarrow_{\frac{1}{n}}B,\Delta^*)[\bar m/\bar x]^{\mathfrak{M}}< \epsilon$ and $(\mathcal{L}^*_{\frac{1}{n}}\mid \Gamma^*\Rightarrow_{\frac{1}{n}}\Delta^*)[\bar m/\bar x]^{\mathfrak{M}}< \epsilon$.

    \item Node $\mathcal{L}^*\mid \Gamma^*\Rightarrow \Delta^*$ in $X_{i+1}$ has its parent node $\mathcal{L}^*\mid \Gamma^*\Rightarrow \Delta^*$ in $X_i$. We simply obtain
    $$(\mathcal{L}^*_{\frac{1}{n}}\mid \Gamma^*\Rightarrow_{\frac{1}{n}}\Delta^*)[\bar m/\bar x]^{\mathfrak{M}}<\max_{G\in X_i}(\mathcal{G}_{\frac{1}{n}}[\bar m/\bar x])^{\mathfrak{M}}<\epsilon$$
\end{enumerate}
Together, we have that $(\mathcal{G}_{\frac{1}{n}}[\bar m/\bar x])^{\mathfrak{M}}<\epsilon$ holds for any $\mathcal{G}\in X_{i+1}$. Therefore, for any $\epsilon>0$ there are $\bar m$ such that
$$\max_{G\in X_i}(\mathcal{G}_{\frac{1}{n}}[\bar m/\bar x])^{\mathfrak{M}}<\epsilon$$
It implies the promised inequality $$\inf_{\bar m\in M}\max_{G\in X_i}(\mathcal{G}_{\frac{1}{n}}[\bar m/\bar x])^{\mathfrak{M}}\leq 0$$

%Left \to
\textbf{Left $\to$ case:}
When the left-most node in $X_{i+1}$ takes the form $\mathcal{L}\mid \Gamma,B\Rightarrow A,\Delta\mid \Gamma\Rightarrow\Delta$ with parent node $\mathcal{L}\mid \Gamma ,A\to B\Rightarrow\Delta$ in $X_i$, the induction hypothesis give
$$\inf_{\bar m}\max_{G\in X_i}(\mathcal{G}_{\frac{1}{n}}[\bar m/\bar x])^{\mathfrak{M}}\leq 0$$
Pick $\bar m$ for each $\epsilon>0$ such that
$$(\mathcal{L}_{\frac{1}{n}}\mid \Gamma\Rightarrow_{\frac{1}{n}} A\to B,\Delta)[\bar m/\bar x]^{\mathfrak{M}}\leq \max_{G\in X_i}(\mathcal{G}_{\frac{1}{n}}[\bar m/\bar x])^{\mathfrak{M}}< \epsilon$$
Since $(A\to B)^{\mathfrak{M}}$ on the left of $\Rightarrow$ and the hyper-connective $\mid$ are both interpreted as minimum, we compute
\begin{eqnarray*}
&&(\mathcal{L}_{\frac{1}{n}}\mid \Gamma,A\Rightarrow_{\frac{1}{n}}B,\Delta\mid \Gamma\Rightarrow_{\frac{1}{n}}\Delta)[\bar m/\bar x]^{\mathfrak{M}}\\
&=&\min((\mathcal{L}_{\frac{1}{n}}\mid \Gamma,A\Rightarrow_{\frac{1}{n}}B,\Delta)[\bar m/\bar x]^{\mathfrak{M}}, (\mathcal{L}_{\frac{1}{n}}\mid \Gamma\Rightarrow_{\frac{1}{n}}\Delta)[\bar m/\bar x]^{\mathfrak{M}})\\
&=&(\mathcal{L}_{\frac{1}{n}}\mid \Gamma,A\to B\Rightarrow_{\frac{1}{n}}\Delta)[\bar m/\bar x]^{\mathfrak{M}}\\
&<& \epsilon
\end{eqnarray*}

For other element $\mathcal{R}$ in $X_{i+1}$, we have $(\mathcal{R}_{\frac{1}{n}}[\bar m/\bar x])^{\mathfrak{M}}<\epsilon$ by a similar argument as before. Thus, we have for each $\epsilon$ there exists $\bar m$ such that $$\max_{G\in X_i}(\mathcal{G}_{\frac{1}{n}}[\bar m/\bar x])^{\mathfrak{M}}<\epsilon$$
and consequently $$\inf_{\bar m\in M}\max_{G\in X_i}(\mathcal{G}_{\frac{1}{n}}[\bar m/\bar x])^{\mathfrak{M}}\leq 0$$

%Right \exists
\textbf{Right $\exists$ case:} When the left-most node in $X_{i+1}$ takes the form $\mathcal{L}\mid \Gamma\Rightarrow A[x/y],\Delta$ with parent node $\mathcal{L}\mid \Gamma \Rightarrow (\exists y)A,\Delta$ in $X_i$, pick $\bar m$ for each $\epsilon>0$ by induction hypothesis such that
$$(\mathcal{L}_{\frac{1}{n}}\mid \Gamma\Rightarrow_{\frac{1}{n}} (\exists y)A,\Delta)[\bar m/\bar x]^{\mathfrak{M}}\leq \max_{G\in X_i}(\mathcal{G}_{\frac{1}{n}}[\bar m/\bar x])^{\mathfrak{M}}< \frac{\epsilon}{2}$$
Since $(\exists y)A$ is interpreted as an infinimum, there exists some $m\in M$ satisfying $$A[m/y,\bar m/\bar x]^{\mathfrak{M}}< (\exists y) A[\bar m/\bar x]^{\mathfrak{M}}+\frac{\epsilon}{2}$$ We then compute
    \begin{eqnarray*} 
        (\mathcal{L}_{\frac{1}{n}}\mid \Gamma\Rightarrow_{\frac{1}{n}}A,\Delta)[m/y,\bar m/\bar x]^{\mathfrak{M}}&<&
        (\mathcal{L}_{\frac{1}{n}}\mid \Gamma\Rightarrow_{\frac{1}{n}}(\exists y)A,\Delta)[\bar m/\bar x]^{\mathfrak{M}}+\frac{\epsilon}{2} \\
        &<& \frac{\epsilon}{2}+\frac{\epsilon }{2}\\
        &=& \epsilon
    \end{eqnarray*}

For other element $\mathcal{R}$ in $X_{i+1}$, there are two possibilities.
\begin{enumerate}
    \item Node $\mathcal{L}^*\mid \Gamma^*\Rightarrow A[x/y],\Delta^*$ in $X_{i+1}$ with its parent node $\mathcal{L}^*\mid \Gamma^*\Rightarrow (\exists y)A,\Delta^*$ in $X_i$. We have by induction hypothesis that
    $$(\mathcal{L}^*_{\frac{1}{n}}\mid \Gamma^*\Rightarrow_{\frac{1}{n}}(\exists y)A,\Delta^*)[\bar m/\bar x]^{\mathfrak{M}}\leq \max_{G\in X_i}(\mathcal{G}_{\frac{1}{n}}[\bar m/\bar x])^{\mathfrak{M}}<\frac{\epsilon}{2}$$
    By using the same $m$, we similarly yield the desired inequality $(\mathcal{L}_{\frac{1}{n}}\mid \Gamma\Rightarrow_{\frac{1}{n}}A[x/y]\Delta)[m/x,\bar m/\bar x]^{\mathfrak{M}}< \epsilon$.

    \item Node $\mathcal{L}^*\mid \Gamma^*\Rightarrow \Delta^*$ in $X_{i+1}$ has its parent node $\mathcal{L}^*\mid \Gamma^*\Rightarrow \Delta^*$ in $X_i$. We simply obtain
    $$(\mathcal{L}^*_{\frac{1}{n}}\mid \Gamma^*\Rightarrow_{\frac{1}{n}}\Delta^*)[\bar m/\bar x]^{\mathfrak{M}}\leq\max_{G\in X_i}(\mathcal{G}_{\frac{1}{n}}[\bar m/\bar x])^{\mathfrak{M}}<\epsilon$$
\end{enumerate}
Together, we conclude that $(\mathcal{G}_{\frac{1}{n}}[m/x,\bar m/\bar x])^{\mathfrak{M}}<\epsilon$ holds for each $\mathcal{G}\in X_{i+1}$. Therefore, for any $\epsilon>0$ there are tuples $(\bar m,m)$ such that
$$\max_{G\in X_i}(\mathcal{G}_{\frac{1}{n}}[m/x,\bar m/\bar x])^{\mathfrak{M}}<\epsilon$$
It implies the promised inequality $$\inf_{m\in M}\inf_{\bar m\in M}\max_{G\in X_i}(\mathcal{G}_{\frac{1}{n}}[\bar m/\bar x])^{\mathfrak{M}}\leq 0$$

%Left \exists

\textbf{Left $\exists$ case:}
When the left-most node in $X_{i+1}$ takes the form $\mathcal{L}\mid \Gamma,A[f(\bar x)/y]\Rightarrow \Delta$ with parent node $\mathcal{L}\mid \Gamma ,(\exists y)A\Rightarrow\Delta$ in $X_i$, the induction hypothesis give
$$\inf_{\bar m}\max_{G\in X_i}(\mathcal{G}_{\frac{1}{n}}[\bar m/\bar x])^{\mathfrak{M}}\leq 0$$
Pick $\bar m$ for each $\epsilon>0$ such that
$$(\mathcal{L}_{\frac{1}{n}}\mid \Gamma,(\exists y)A\Rightarrow_{\frac{1}{n}}\Delta)[\bar m/\bar x]^{\mathfrak{M}}\leq \max_{G\in X_i}(\mathcal{G}_{\frac{1}{n}}[\bar m/\bar x])^{\mathfrak{M}}< \epsilon$$
Since $(\exists y)A$ is interpreted as a supremum on the left of $\Rightarrow$, we have that 
$$(\mathcal{G}_{\frac{1}{n}}\mid \Gamma, A[f(\bar x)/y]\Rightarrow_{\frac{1}{n}}\Delta)[\bar m/\bar x]^{\mathfrak{M}}\leq (\mathcal{G}_{\frac{1}{n}}\mid \Gamma, (\exists y)A\Rightarrow_{\frac{1}{n}}\Delta)[\bar m/\bar x]^{\mathfrak{M}}< \epsilon$$ 

For other element $\mathcal{R}$ in $X_{i+1}$, we have $(\mathcal{R}_{\frac{1}{n}}[\bar m/\bar x])^{\mathfrak{M}}<\epsilon$ by a similar argument. Thus, we have for each $\epsilon$ there exists $\bar m$ such that $$\max_{G\in X_i}(\mathcal{G}_{\frac{1}{n}}[\bar m/\bar x])^{\mathfrak{M}}<\epsilon$$
and consequently $$\inf_{\bar m\in M}\max_{G\in X_i}(\mathcal{G}_{\frac{1}{n}}[\bar m/\bar x])^{\mathfrak{M}}\leq 0$$
\end{proof}
\end{lemma}

Notice that $X_{h}$, the leaf set of $T_{\mathcal{H}}$, consists of hypersequents with no quantifiers. For any given $n$, the above Lemma \ref{apsa} applying to $X_{h}$ yields
$$\sup_{\mathfrak{M}}\inf_{\bar m\in M} \max_{\mathcal{G}\in X_{h}}(\mathcal{G}_{\frac{1}{n}}[\bar m/\bar x])^{\mathfrak{M}}\leq 0$$ 
It can be viewed as a $\Sigma_1^0$-judgment stating $$``(\exists \bar x) (\bigwedge_{\mathcal{G}\in X_{h}} \mathcal{G}_{\frac{1}{n}})"$$ 
We will adopt an approximate version of Herbrand's theorem to show existence of Herbrand witnesses. 
  
\begin{theorem}[Approximate Herbrand Theorem]\label{aphb}
Suppose $X$ is a finite set of quantifier-free hypersequents with shared free variables $\bar x$ and $n$ is a natural number. Then there are terms $\bar t_1,\ldots, \bar t_k$ such that 

$$\sup_{\mathfrak{M}}\min_{1\leq i\leq k} \max_{\mathcal{G}\in X}(\mathcal{G}[\bar t_i/\bar x]^{\mathfrak{M}})< \sup_{\mathfrak{M}}\inf_{{\bar m}\in M} \max_{\mathcal{G}\in X}(\mathcal{G}[\bar m/\bar x]^{\mathfrak{M}})+\frac{1}{n}$$

If we additionally assume $$\sup_{\mathfrak{M}}\inf_{{\bar m}\in M} \max_{\mathcal{G}\in X}(\mathcal{G}[\bar m/\bar x]^{\mathfrak{M}})\leq 0$$ then the hypersequent $$\mathcal{G}^1_{\frac{1}{n}}[\bar t_1/\bar x]\mid \cdots \mid \mathcal{G}^k_{\frac{1}{n}}[\bar t_k/\bar x]$$ is valid for each $\mathcal{G}^1,\ldots,\mathcal{G}^k\in X$. 
\begin{proof}
We generalize the proof appearing in \cite{main}. Let $P$ collect all formulas of the form $r(\bar t)$ where $r$ is a relation symbol and $\bar t$ are terms in the language. Note that a structure $\mathfrak{M}$ is uniquely associated with a map $v_{\mathfrak{M}}\in [0,1]^{P}$. For simplicity, let
$$d:=\sup_{\mathfrak{M}}\inf_{{\bar m}\in M} \max_{\mathcal{G}\in X}(\mathcal{G}[\bar m/\bar x])^{\mathfrak{M}}$$

Given arbitrary term $\bar t$, we define 
\begin{eqnarray*}
S(\bar t)&:=&\{v_{\mathfrak{M}}\in [0,1]^{P}: \max_{\mathcal{G}\in X}(\mathcal{G}[\bar t/ \bar x]^{\mathfrak{M}})\geq d+\frac{1}{n}\}\\
&=&\bigcup_{\mathcal{G}\in X}\{v_{\mathfrak{M}}\in [0,1]^{P}: (\mathcal{G}[\bar t/ \bar x]^{\mathfrak{M}})\geq d+\frac{1}{n}\}
\end{eqnarray*}

All $S(\bar t)$ are closed subsets of $[0,1]^{P}$ because $\to$, $\Rightarrow$, and $\mid $ are interpreted as continuous functions on $\mathbb{R}$ with absence of quantifiers in $\mathcal{G}$ and finite unions of closed sets are still closed. Define $$\mathcal{S}:=\{S(\vec t)\mid \bar t \text{ is a term}\}\subseteq \mathcal{P}([0,1]^{P})$$

Suppose towards a contradiction that there are no finite list of terms $\bar t_1,\ldots,\bar t_k$ satisfying $\sup_{\mathfrak{M}}\min_i \max_{\mathcal{G}\in X}  (\mathcal{G}[\bar t_i/\bar x])^\mathfrak{M}<d+\frac{1}{n}$. Then for any given list of terms $\bar t_1,\ldots,\bar t_k$, there exists a structure $\mathfrak{M}$ such that $$\max_{\mathcal{G}\in X}(\mathcal{G}[\bar t_i/\bar x])^\mathfrak{M}\geq d+\frac{1}{n}$$ for all $1\leq i\leq k$. Hence, $v_{\mathfrak{M}}\in \bigcap_i S(\bar t_i)$ and we have shown that $\mathcal{S}$ has the finite intersection property. 

Since $\mathcal{S}$ is a collection of closed subsets of $[0,1]^{P}$ with finite intersection property and $[0,1]^{P}$ is a compact space by Tychonoff's Theorem, $\bigcap \mathcal{S}$ is nonempty and thus there exists a structure $\mathfrak{M}$ such that $\max_{\mathcal{G}\in X}(\mathcal{G}[\bar t/\bar x]^{\mathfrak{M}})\geq d+\frac{1}{n}$ for all terms $\bar t$. It would particularly imply that 
\begin{eqnarray*}
\inf_{\bar m\in M}\max_{\mathcal{G}\in X}(\mathcal{G}[\bar m/\bar x]^{\mathfrak{M}})&\geq& d+\frac{1}{n}\\
&=& \sup_{\mathfrak{M}}\inf_{{\bar m}\in M} \max_{\mathcal{G}\in X}(\mathcal{G}[\bar m/\bar x])^{\mathfrak{M}}+\frac{1}{n}\\
&>&\sup_{\mathfrak{M}}\inf_{{\bar m}\in M} \max_{\mathcal{G}\in X}(\mathcal{G}[\bar m/\bar x])^{\mathfrak{M}}
\end{eqnarray*}
which itself is a contradiction. Therefore, there must exist a list of Herbrand terms $\bar t_1,\ldots, \bar t_k$ as promised. 

If furthermore we have $\sup_{\mathfrak{M}}\inf_{{\bar m}\in M} \max_{\mathcal{G}\in X}(\mathcal{G}[\bar m/\bar x])^{\mathfrak{M}}\leq 0$, then for any given $\mathcal{G}^1,\ldots,\mathcal{G}^k$ in $X$ we have that
$$(\mathcal{G}^i_{\frac{1}{n}}[\bar t_i/\bar x]^{\mathfrak{M}}) 
\leq \max_{\mathcal{G}\in X} (\mathcal{G}_{\frac{1}{n}}[\bar t_i/\bar x]^{\mathfrak{M}} $$
Therefore, 
\begin{eqnarray*}
\sup_{\mathfrak{M}}\min_{1\leq i\leq k}(\mathcal{G}^i_{\frac{1}{n}}[\bar t_i/\bar x]^{\mathfrak{M}}) 
&\leq& \sup_{\mathfrak{M}}\min_{i}\max_{\mathcal{G}\in X} (\mathcal{G}_{\frac{1}{n}}[\bar t_i/\bar x]^{\mathfrak{M}}) \\
&=&\sup_{\mathfrak{M}}\min_{i} \max_{\mathcal{G}\in X}(n\cdot (\mathcal{G}[\bar t_i/ \bar x]^{\mathfrak{M}})-1)\\
&=& n\cdot  \sup_{\mathfrak{M}}\min_{i} \max_{\mathcal{G}\in X}( (\mathcal{G}[\bar t_i/ \bar x]^{\mathfrak{M}})-\frac{1}{n})\\
&<& n\cdot\sup_{\mathfrak{M}}\inf_{{\bar m}\in M} \max_{\mathcal{G}\in X}(\mathcal{G}[\bar m/\bar x])^{\mathfrak{M}}\\
&\leq& 0
\end{eqnarray*}
which exactly means that $\mathcal{G}^1_{\frac{1}{n}}[\bar t_1/\bar x]\mid \cdots \mid \mathcal{G}^k_{\frac{1}{n}}[\bar t_k/\bar x]$ is a valid hypersequent.
\end{proof}
\end{theorem}

Applying the above theorem to $X_{h}$, we arrive at some Herbrand witnesses $\bar t_1,\ldots, \bar t_k$ such that the following hypersequent consisting of only atomic formulas $$\mathcal{G}^1_{\frac{1}{n}}[\bar t_1/\bar x]\mid \cdots \mid \mathcal{G}^k_{\frac{1}{n}}[\bar t_k/\bar x]$$
is valid for each $\mathcal{G}^1,
\ldots,\mathcal{G}^k$ in $X_{h}$. From now on, we omit $\bar x$ and denote it as $$\mathcal{G}^1_{\frac{1}{n}}(\bar t_1)\mid \cdots \mid \mathcal{G}^k_{\frac{1}{n}}(\bar t_k)$$ to simplify the notation.

These are derivable hypersequents by completeness of GŁ. We will recover $\mathcal{H}_{\frac{1}{n}}$ by re-introducing connectives following the reverse order of the Skolemization tree $T_{\mathcal{H}}$. For that purpose, we need the following admissible rules that handle introduction rules for connectives in the sense of $\Rightarrow_{\frac{1}{n}}$.

\begin{lemma}
The followings are derivable rules in GŁ and GŁ$\forall$ for any $n$.

\begin{center}
    \AxiomC{$\mathcal{G}\mid \Gamma\Rightarrow \Delta$}
     \RightLabel{mul$\qquad\qquad \qquad$}   
    \UnaryInfC{$\mathcal{G}\mid n\Gamma\Rightarrow n\Delta $}
\DisplayProof
    \AxiomC{$\mathcal{G}\mid n\Gamma\Rightarrow n\Delta$}
     \RightLabel{div}   
    \UnaryInfC{$\mathcal{G}\mid \Gamma\Rightarrow \Delta $}
\DisplayProof
\end{center}
\begin{proof}
The $mul$ rule is simply an application of $mix$ of $n$-copies of $\mathcal{G}\mid \Gamma \Rightarrow \Delta$. The $div$ is the result of applying $split$ and $ec$. 
\end{proof}
\end{lemma}

\begin{lemma}\label{to1n}
We have left and right introduction rules for $\to$ in the sense of $\Rightarrow_\frac{1}{n}$ for any $n$.

\begin{center}
    \AxiomC{$\mathcal{G}\mid \Gamma,B\Rightarrow_{\frac{1}{n}}A,\Delta\mid \Gamma \Rightarrow_{\frac{1}{n}}\Delta$}
    \RightLabel{$\to \Rightarrow_{\frac{1}{n}}$}
    \UnaryInfC{$\mathcal{G}\mid \Gamma,A\to B\Rightarrow_{\frac{1}{n}}\Delta$}
\DisplayProof
    \AxiomC{$\mathcal{G}\mid \Gamma,A\Rightarrow_{\frac{1}{n}}B,\Delta$}
    \AxiomC{$\mathcal{G}\mid \Gamma\Rightarrow_{\frac{1}{n}}\Delta$}
     \RightLabel{$ \Rightarrow_{\frac{1}{n}}\to $}
    \BinaryInfC{$\mathcal{G}\mid \Gamma\Rightarrow_{\frac{1}{n}}A\to B,\Delta$}
\DisplayProof
\end{center}

\begin{proof}
The $\to\Rightarrow_{\frac{1}{n}}$ rule can be derived by iterating the following derivable rules and apply $ec$.

\begin{prooftree}
    \AxiomC{$\mathcal{G}\mid \Gamma,B\Rightarrow A,\Delta$}
    \RightLabel{$ew$}
    \UnaryInfC{$\mathcal{G}\mid\Gamma,B\Rightarrow A,\Delta\mid \Gamma\Rightarrow\Delta$}
        \RightLabel{$\to\Rightarrow$}
    \UnaryInfC{$\mathcal{G}\mid\Gamma,A\to B\Rightarrow \Delta$}
\end{prooftree}
and 

\begin{prooftree}
    \AxiomC{$\mathcal{G}\mid \Gamma\Rightarrow \Delta$}
    \RightLabel{$wl$}
    \UnaryInfC{$\mathcal{G}\mid\Gamma,A\to B\Rightarrow \Delta$}
\end{prooftree}

For the $\Rightarrow_{\frac{1}{n}}\to$ rule, we first notice that $\mathcal{G}\mid \bot, n\Gamma,iA\Rightarrow iB,n\Delta$ is derivable for any $0<i<n$.
\begin{prooftree}
    \AxiomC{$\mathcal{G}\mid \bot, n\Gamma,nA\Rightarrow nB,n\Delta$}
    \RightLabel{$mul$}
    \UnaryInfC{$\mathcal{G}\mid i\bot, ni\Gamma,ni A\Rightarrow ni B,nk\Delta$}
    \AxiomC{$\mathcal{G}\mid \bot, n\Gamma\Rightarrow n\Delta$}
    \RightLabel{$mul$}
     \UnaryInfC{$\mathcal{G}\mid (n-i)\bot, n(n-i)\Gamma\Rightarrow n(n-i)\Delta$}
     \RightLabel{$mix$}
    \BinaryInfC{$\mathcal{G}\mid n\bot, n^2\Gamma,ni A\Rightarrow niB,n^2\Delta$}
    \RightLabel{$div$}
    \UnaryInfC{$\mathcal{G}\mid \bot, n\Gamma,i A\Rightarrow iB,n\Delta$}
\end{prooftree}
For each $i$ and $j$, we have 
\begin{prooftree}
    \AxiomC{$\mathcal{G}\mid \bot, n\Gamma,(i+1)A\Rightarrow (i+1)B,j(A\to B),n\Delta$}
    \AxiomC{$\mathcal{G}\mid \bot, n\Gamma,iA\Rightarrow iB, j(A\to B), n\Delta$}
    \BinaryInfC{$\mathcal{G}\mid \bot, n\Gamma,i A\Rightarrow iB, (j+1)(A\to B),n\Delta$}
\end{prooftree}
We thus have a derivation of $\mathcal{G}\mid \bot, n\Gamma\Rightarrow  n(A\to B),n\Delta$ by a recursion on $j$.
\end{proof}
\end{lemma}

\begin{lemma}\label{ex1n}
We have left and right introduction rules for $\exists$ in the sense of $\Rightarrow_\frac{1}{n}$ for any $n$. For the left rule, we require that $c$ is not free in $\Gamma$, $\Delta$, and $\mathcal{G}$.

\begin{center}
    \AxiomC{$\mathcal{G}\mid \Gamma,A(c)\Rightarrow_{\frac{1}{n}} \Delta$}
    \RightLabel{$\exists\Rightarrow_{\frac{1}{n}}\qquad$}
    \UnaryInfC{$\mathcal{G}\mid \Gamma,(\exists x)A[x/c]\Rightarrow_{\frac{1}{n}} 
    \Delta$}
\DisplayProof
    \AxiomC{$\mathcal{G}\mid \Gamma\Rightarrow_{\frac{1}{n}} A(t),\Delta$}
     \RightLabel{$\Rightarrow_{\frac{1}{n}}\exists$}
    \UnaryInfC{$\mathcal{G}\mid \Gamma\Rightarrow_{\frac{1}{n}}(\exists x)A[x/t], \Delta$}
\DisplayProof
\end{center}

(This lemma and a similar proof also appeared in \cite{gera}).

\begin{proof}
In the case of $\exists\Rightarrow_{\frac{1}{n}}$, take fresh new variables $x_1,\ldots, x_n$. We have deductions $\mathcal{G}\mid \bot, n\Gamma,n A(x_i)\Rightarrow n\Delta$ for each $i$ by substituting $x$ with $x_i$. We thus have the following derivation for $\mathcal{G}\mid \bot, n\Gamma,n(\exists x)A\Rightarrow n\Delta$. 

\begin{prooftree}
    \AxiomC{$\overbrace{\mathcal{G}\mid \bot, n\Gamma,n A (x_1)\Rightarrow n\Delta \qquad \cdots\qquad\mathcal{G}\mid \bot, n\Gamma,n A(x_n)\Rightarrow n\Delta}^n$}
    \RightLabel{$mix$}
    \UnaryInfC{$\mathcal{G}\mid n\bot, n^2\Gamma,n^2A(x_1),\ldots,n^2A(x_n)\Rightarrow n^2\Delta$}
    \RightLabel{$split$}
    \UnaryInfC{$\mathcal{G}\mid \overbrace{\bot, n\Gamma,A(x_1),\ldots,A(x_n)\Rightarrow n\Delta\mid \cdots \mid \bot, n\Gamma,A(x_1),\ldots,A(x_n)\Rightarrow n\Delta}^{n}$}
    \RightLabel{$ec$}
    \UnaryInfC{$\mathcal{G}\mid \bot, n\Gamma,A(x_1),\ldots,A(x_n)\Rightarrow n\Delta$}
    \RightLabel{$\exists\Rightarrow$}
    \UnaryInfC{$\mathcal{G}\mid \bot, n\Gamma,n(\exists x)A\Rightarrow n\Delta$}
\end{prooftree}

The $\Rightarrow_{\frac{1}{n}}\exists$ rule can be derived by simply iterating $\Rightarrow\exists$, since there is no eigen-conditions. 
\end{proof}
\end{lemma}

Now we are ready to present the main completeness result. We need a few definitions on multi-labeled hypersequents where labels keep trace of the order of connective introduction.

\begin{definition}
Consider a valid hypersequent $\mathcal{H}$ with no free variables, its tree $T_{\mathcal{H}}$ and leafs $X_{h}$ defined as above, Herbrand witnesses $\bar t_1,\ldots, \bar t_k$ for $\bar x$ from Lemma \ref{aphb}, and arbitrary natural number $n$.

A \textbf{labeled hypersequent} is of the form $\mathcal{G}^{\sigma}$ where $\mathcal{G}$ is a hypersequent and $\sigma$ is a path in $T_{\mathcal{H}}$ such that the following holds
\begin{itemize}
    \item $\mathcal{G} = \mathcal{L}_{\frac{1}{n}}(\bar t_i)$ where $\mathcal{L}$ is the last node of $\sigma$ and $\bar t_i$ is some Herbrand witness.
\end{itemize}
A \textbf{multi-labeled hypersequent} $\mathfrak{G}$ is a multiset of labeled hypersequents of the form $$\mathcal{G}_1^{\sigma_1}\mid \cdots\mid \mathcal{G}_l^{\sigma_l}$$
where $0\leq l\leq k$. We also refer $\mathcal{G}_i^{\sigma_i}$ as a \textbf{component} of $\mathfrak{G}$.

We say a labeled hypersequent $\mathcal{G}^{\sigma}$ is \textbf{derivable} if $\mathcal{G}$ is derivable; a multi-labeled hypersequent $\mathcal{G}_1^{\sigma_1}\mid \cdots\mid \mathcal{G}_k^{\sigma_k}$ is \textbf{derivable} if $\mathcal{G}_1\mid \cdots\mid \mathcal{G}_k$ is derivable. 

The \textbf{stage} of a multi-labeled hypersequent $\mathcal{G}_1^{\sigma_1}\mid \cdots\mid \mathcal{G}_k^{\sigma_k}$ is the tuple $(|\sigma_1|,\ldots,|\sigma_k|)$ where $|\sigma_i|$ is the length of path $\sigma_i$ in tree $T_{\mathcal{H}}$. The \textbf{size} of $\mathcal{G}_1^{\sigma_1}\mid \cdots\mid \mathcal{G}_k^{\sigma_k}$ is $\sum_{1\leq i\leq k} |\sigma_i|$.

\end{definition}

\begin{definition}\label{wlfm}
Given a finite set $\mathbb{X}$ of multi-labeled hypersequents, we say $\mathbb{X}$ is \textbf{well-formed} if the following conditions hold:
\begin{enumerate}
    \item There exists some element $\mathfrak{G}$ such that each label in $\mathfrak{G}$ is an initial segment of the left-most path in $T_{\mathcal{H}}$. We call $\mathfrak{G}$ the \textbf{left-most element}.
    
    \item For any element $\mathfrak{G}\mid \mathcal{L}^\sigma\mid\mathfrak{L}$ in $\mathbb{X}$ where $\sigma$ extends some path $\eta \mathbin{\frown} (\mathcal{G}\mid \Gamma\Rightarrow A\to B, \Delta)\mathbin{\frown}(\mathcal{G} \mid \Gamma,A \Rightarrow B,\Delta)$ in $T_{\mathcal{H}}$, there exists some element $\mathfrak{G}\mid \mathcal{R}^{\tau}\mid \mathfrak{L}$ in $\mathbb{X}$ such that $\tau$ extends
$\eta \mathbin{\frown} (\mathcal{G}\mid \Gamma\Rightarrow A\to B, \Delta)\mathbin{\frown}(\mathcal{G} \mid \Gamma \Rightarrow \Delta)$.

%     \item For any element $\mathfrak{G}\mid \mathcal{R}^\tau$ in $\mathbb{X}$ where $\tau$ extends some path $\eta \mathbin{\frown} (\mathcal{G}\mid \Gamma\Rightarrow A\to B, \Delta)\mathbin{\frown}(\mathcal{G} \mid \Gamma \Rightarrow \Delta)$ in $T$, there exists some element $\mathfrak{G}\mid \mathcal{L}^{\sigma}$ in $\mathbb{X}$ such that $\sigma$ extends
% $\eta \mathbin{\frown} (\mathcal{G}\mid \Gamma\Rightarrow A\to B, \Delta)\mathbin{\frown}(\mathcal{G} \mid \Gamma ,A\Rightarrow B,\Delta)$.
    
    \item Every element of $\mathbb{X}$ has the same number of components and the same stage. 

    \item Every element of $\mathbb{X}$ is derivable in GŁ$\forall$.
\end{enumerate}
We define $\mu(\mathbb{X})$ the \textbf{size} of $\mathbb{X}$ to be the size of its left-most element.
\end{definition}

The prototype of a well-formed set is the following:
\begin{definition}\label{init}
Define the following set of multi-labeled hypersequents
$$\mathbb{X}_0:=\{(\mathcal{G}^1_{\frac{1}{n}}(\bar t_1))^{\sigma_1}\mid\cdots  \mid (\mathcal{G}^k_{\frac{1}{n}}(\bar t_k))^{\sigma_k}: \mathcal{G}^i\in X_{h}\text{ and }\sigma_i \text{ is the path ending with }\mathcal{G}^i\}$$
and a left-most element in $\mathbb{X}_0$ $$\mathfrak{G}_0:=(\mathcal{G}_{\frac{1}{n}}(\bar t_1))^{\sigma}\mid\cdots  \mid (\mathcal{G}_{\frac{1}{n}}(\bar t_k))^{\sigma}$$
where $\sigma$ is the left-most path in $T_{\mathcal{H}}$ ending with $\mathcal{G}$. 
\end{definition}

\begin{lemma}\label{inwf}
The set $\mathbb{X}_0$ is well-formed.
\begin{proof}
We verify each condition in Definition \ref{wlfm}.
\begin{enumerate}
    \item $\mathfrak{G}_0$ is by definition a left-most element.
    \item $\mathbb{X}_0$ enumerates all paths in $T_{\mathcal{H}}$ with any combination of witnesses. Hence, for any $\mathfrak{G}\mid \mathcal{G}^{\sigma}$ in $\mathbb{X}_0$ where $\sigma$ extends some left branch at a split in $T_{\mathcal{H}}$, say $\tau$ is a path in $T_{\mathcal{H}}$ extending the right branch and ending with $\mathcal{R}$. Then $\mathfrak{G}\mid \mathfrak{L}^{\tau}$ is also in $\mathbb{X}_0$ by design.
    \item Elements of $\mathbb{X}_0$ all have the same stage $(h,\ldots,h)$ by Lemma \ref{sync} and $k$-many components.
    \item Since $\mathbb{X}_0$ consists of valid quantifier-free hypersequents by Lemma \ref{aphb}, its elements are derivable in GŁ by Theorem \ref{GŁcmp} and thus also derivable in GŁ$\forall$.
\end{enumerate}
\end{proof}
\end{lemma}

Starting from $\mathbb{X}_0$, we will iteratively perform a reduction of labels in $\mathbb{X}_0$ until all labels are of the form $(\mathcal{H})$.

\begin{definition}\label{lbrd}
We define a function $F$ on well-formed sets of multi-labeled hypersequents $\mathbb{X}$ with left-most element $\mathfrak{G}$ by the following. The order can be arbitrarily chosen if multiple conditions are satisfied.

% Left \to
\textbf{(1):} If $\mathfrak{G}$ contains a component of the form $\mathcal{L}^{\sigma}$ where 
    \begin{eqnarray*}
     &&\sigma =\eta \mathbin{\frown} (\mathcal{G}\mid \Gamma, A\to B\Rightarrow \Delta)\mathbin{\frown}(\mathcal{G} \mid \Gamma,B \Rightarrow A,\Delta\mid \Gamma \Rightarrow \Delta)\\
    &&\mathcal{L} =(\mathcal{G}_{\frac{1}{n}} \mid \Gamma,B \Rightarrow_{\frac{1}{n}} A,\Delta\mid \Gamma \Rightarrow_{\frac{1}{n}} \Delta)(\bar t_i)
    \end{eqnarray*}
For simplicity, we write $\mathfrak{G}$ as $\mathfrak{L}\mid \mathcal{L}^{\sigma}$ by permuting $\mathcal{L}^{\sigma}$ to the right since multi-labeled hypersequents are multisets of labeled hypersequents.

By well-formedness of $\mathbb{X}$, $\mathfrak{L}\mid\mathcal{L}^{\sigma}$ is derivable. Now consider the following proof using Lemma \ref{to1n} (without labels)
\begin{prooftree}
    \AxiomC{$\mathfrak{L}\mid (\mathcal{G}_{\frac{1}{n}} \mid \Gamma,B \Rightarrow_{\frac{1}{n}} A,\Delta\mid \Gamma \Rightarrow_{\frac{1}{n}} \Delta)(\bar t_i)$}
    \RightLabel{$\to \Rightarrow_{\frac{1}{n}}$}
    \UnaryInfC{$\mathfrak{L}\mid (\mathcal{G}_{\frac{1}{n}} \mid \Gamma,A\to B \Rightarrow_{\frac{1}{n}} \Delta)(\bar t_i)$}
\end{prooftree}
We first replace $\mathfrak{G}$ in $\mathbb{X}$ with the derived hypersequent with a shorter label
$$\mathfrak{L}\mid ((\mathcal{G}_{\frac{1}{n}} \mid \Gamma,A\to B \Rightarrow_{\frac{1}{n}} \Delta)(\bar t_i))^{\eta \mathbin{\frown} (\mathcal{G}\mid \Gamma, A\to B\Rightarrow \Delta)}$$

For other elements $\mathfrak{R}$ in $\mathbb{X}$ with component $(\mathcal{L}^*)^{\sigma^*}$ corresponding to $\mathcal{L}^{\sigma}$ in $\mathfrak{G}$, we denote $\mathfrak{R}$ as $\mathfrak{L}^*\mid (\mathcal{L}^*)^{\sigma^*}$. There are two cases due to the construction in $T_{\mathcal{H}}$ in Definition \ref{dctr}. 
        \begin{itemize}
    \item When $ (\mathcal{L}^*)^{\sigma^*}$ takes the following form 
    \begin{eqnarray*}
     &&\sigma^* =\eta^* \mathbin{\frown} (\mathcal{G}^*\mid \Gamma^*, A\to B\Rightarrow \Delta^*)\mathbin{\frown}(\mathcal{G}^* \mid \Gamma^*,B \Rightarrow A,\Delta^*\mid \Gamma^* \Rightarrow \Delta^*)\\
    &&\mathcal{L}^* = (\mathcal{G}^*_{\frac{1}{n}} \mid \Gamma^*,B \Rightarrow_{\frac{1}{n}} A,\Delta^*\mid \Gamma^* \Rightarrow_{\frac{1}{n}} \Delta^*)(\bar t_i)
    \end{eqnarray*}
We derive the following from $\mathfrak{L}^*\mid \mathcal{L}^*$ similarly to the case for $\mathfrak{G}$
$$\mathfrak{L}^*\mid ((\mathcal{G}^*_{\frac{1}{n}} \mid \Gamma^*,A\to B \Rightarrow_{\frac{1}{n}} \Delta^*)(\bar t_i))^{\eta^* \mathbin{\frown} (\mathcal{G}^*\mid \Gamma^*, A\to B\Rightarrow \Delta^*)}$$
and replace $\mathfrak{R}$ in $\mathbb{X}$ with it.

    \item When $ (\mathcal{L}^*)^{\sigma^*}$ takes the following form
    $\mathfrak{L}^*\mid (\mathcal{L}^*)^{\sigma^*}$ where 
\begin{eqnarray*}
    &&\sigma^* =\eta^* \mathbin{\frown} (\mathcal{G}^*\mid \Gamma^*\Rightarrow \Delta^*)\mathbin{\frown}(\mathcal{G}^* \mid \Gamma^*\Rightarrow \Delta^*)\\
   &&\mathcal{L}^* = (\mathcal{G}^*_{\frac{1}{n}} \mid \Gamma^* \Rightarrow_{\frac{1}{n}} \Delta^*)(\bar t_i)
    \end{eqnarray*}
We replace $\mathcal{R}$ with
$$\mathfrak{L}^*\mid ((\mathcal{G}^*_{\frac{1}{n}} \mid \Gamma^* \Rightarrow_{\frac{1}{n}} \Delta^*)(\bar t_i))^{\eta^* \mathbin{\frown} (\mathcal{G}^*\mid \Gamma^*\Rightarrow \Delta^*)}$$
        \end{itemize}
        
Now we define $F(\mathbb{X})$ to be the result of substituting $\mathfrak{G}$ and other elements in $\mathbb{X}$. We check $F(\mathbb{X})$ is still well-formed. 
\begin{itemize}
    \item Property (1) is clear as the replacement of $\mathfrak{G}$ is still a left-most element. 
    \item The only way Property (2) fails is when we remove labels on the right branch of a split while keeping the left branch. Since we do not remove any label of a split at this stage, Property (2) holds for $F(\mathbb{X})$ given that it holds for $\mathbb{X}$. 
    \item Elements in $F(\mathbb{X})$ are kept with the same number of components and the same stage. 
    \item Elements in $F(\mathbb{X})$ are all provided with a derivation.
\end{itemize}

% Right \exists  
\textbf{(2):}
If the condition in $(1)$ fails and $\mathfrak{G}$ is of the form $\mathfrak{L}\mid \mathcal{L}^\sigma$ where 
    \begin{eqnarray*}
     &&\sigma =\eta \mathbin{\frown} (\mathcal{G}\mid \Gamma\Rightarrow \exists y A, \Delta)\mathbin{\frown}(\mathcal{G} \mid \Gamma \Rightarrow A[x/y],\Delta)\\
    &&\mathcal{L}= (\mathcal{G}_{\frac{1}{n}} \mid \Gamma \Rightarrow_{\frac{1}{n}} A[x/y],\Delta)(\bar t_i)
    \end{eqnarray*}
By well-formedness of $\mathbb{X}$, $\mathfrak{L}\mid\mathcal{L}^{\sigma}$ is derivable. Consider the following proof using Lemma \ref{ex1n}
\begin{prooftree}
    \AxiomC{$\mathfrak{L}\mid (\mathcal{G}_{\frac{1}{n}} \mid \Gamma\Rightarrow_{\frac{1}{n}} A[x/y],\Delta)(\bar t_i)$}
    \RightLabel{$ \Rightarrow_{\frac{1}{n}}\exists$}
    \UnaryInfC{$\mathfrak{L}\mid (\mathcal{G}_{\frac{1}{n}} \mid \Gamma \Rightarrow_{\frac{1}{n}} \exists y A,\Delta)(\bar t_i)$}
\end{prooftree}
We then replace $\mathfrak{G}$ in $\mathbb{X}$ with the derived hypersequent with a shorter label
$$\mathfrak{L}\mid ((\mathcal{G}_{\frac{1}{n}} \mid \Gamma\Rightarrow_{\frac{1}{n}}\exists y A, \Delta)(\bar t_i))^{\eta \mathbin{\frown} (\mathcal{G}\mid \Gamma\Rightarrow \exists y A,\Delta)}$$

For other elements $\mathfrak{R}$ in $\mathbb{X}$, we similarly replace it with either 
$$\mathfrak{L}^*\mid ((\mathcal{G}^*_{\frac{1}{n}} \mid \Gamma^*\Rightarrow_{\frac{1}{n}} \exists y A,\Delta^*)(\bar t_i))^{\eta^* \mathbin{\frown} (\mathcal{G}^*\mid \Gamma^*\Rightarrow \exists y A, \Delta^*)}$$
or $$\mathfrak{L}^*\mid ((\mathcal{G}^*_{\frac{1}{n}} \mid \Gamma^*\Rightarrow_{\frac{1}{n}} \Delta^*)(\bar t_i))^{\eta^* \mathbin{\frown} (\mathcal{G}^*\mid \Gamma^*\Rightarrow \Delta^*)}$$
according to its cases. 

Define $F(\mathbb{X})$ to be the result of substituting $\mathfrak{G}$ and other elements in $\mathbb{X}$. By a similar argument as in $(1)$, $F(\mathbb{X})$ remains to be well-formed. 

% Right \to
\textbf{(3):} If conditions in $(1)$ and $(2)$ fail and $\mathfrak{G}$ is of the form $\mathfrak{L}\mid \mathcal{L}^\sigma$ where 
    \begin{eqnarray*}
     &&\sigma =\eta \mathbin{\frown} (\mathcal{G}\mid \Gamma\Rightarrow A\to B, \Delta)\mathbin{\frown}(\mathcal{G} \mid \Gamma,A \Rightarrow B,\Delta)\\
    &&\mathcal{L} = (\mathcal{G}_{\frac{1}{n}} \mid \Gamma ,A\Rightarrow_{\frac{1}{n}} B,\Delta)(\bar t_i)
    \end{eqnarray*}
Notice that $\sigma$ cannot be the right branch of a split due to the choice of $\mathfrak{G}$ being the left-most element. Since $\mathbb{X}$ is well-formed, $\mathfrak{L}\mid\mathcal{L}^{\sigma}$ is derivable. Furthermore, the multi-labeled hypersequent $\mathfrak{L}\mid \mathcal{R}^{\tau}$ where 
    \begin{eqnarray*}
    &&\tau =\eta \mathbin{\frown} (\mathcal{G}\mid \Gamma\Rightarrow A\to B, \Delta)\mathbin{\frown}(\mathcal{G} \mid \Gamma\Rightarrow \Delta)\\
   &&\mathcal{R} = (\mathcal{G}_{\frac{1}{n}} \mid \Gamma \Rightarrow_{\frac{1}{n}} \Delta)(\bar t_i)
    \end{eqnarray*}
exists in $\mathbb{X}$ by condition (2) and (3) of $\mathbb{X}$ being well-formed. Therefore, $\mathfrak{L}\mid \mathcal{R}^{\tau}$ is also derivable. Now consider the following proof using Lemma \ref{to1n}
\begin{prooftree}
    \AxiomC{$\mathfrak{L}\mid (\mathcal{G}_{\frac{1}{n}} \mid \Gamma,A\Rightarrow_{\frac{1}{n}} B,\Delta)(\bar t_i)$}
    \AxiomC{$\mathfrak{L}\mid (\mathcal{G}_{\frac{1}{n}} \mid \Gamma\Rightarrow_{\frac{1}{n}}\Delta)(\bar t_i)$}
    \RightLabel{$ \Rightarrow_{\frac{1}{n}} \to$}
    \BinaryInfC{$\mathfrak{L}\mid (\mathcal{G}_{\frac{1}{n}} \mid \Gamma \Rightarrow_{\frac{1}{n}} A\to B,\Delta)(\bar t_i)$}
\end{prooftree}
We then replace both $\mathfrak{L}\mid \mathcal{L}^{\sigma}$ and $\mathfrak{L}\mid \mathcal{R}^{\tau}$ in $\mathbb{X}$ with 
$$\mathfrak{L}\mid ((\mathcal{G}_{\frac{1}{n}} \mid \Gamma\Rightarrow_{\frac{1}{n}}A\to B, \Delta)(\bar t_i))^{\eta \mathbin{\frown} (\mathcal{G}\mid \Gamma\Rightarrow A\to B,\Delta)}$$

For other elements $\mathfrak{R}$ in $\mathbb{X}$ with $(\mathcal{L}^*)^{\sigma^*}$ corresponding to $\mathcal{L}^{\sigma}$, there are two cases.
\begin{itemize}

    \item $\mathfrak{R}$ is of the form  
    $\mathfrak{L}^*\mid (\mathcal{L}^*)^{\sigma^*}$ where 
    \begin{eqnarray*}
     &&\sigma^* =\eta^* \mathbin{\frown} (\mathcal{G}^*\mid \Gamma^*\Rightarrow A\to B,\Delta^*)\mathbin{\frown}(\mathcal{G}^* \mid \Gamma^*,A \Rightarrow B,\Delta^*)\\
    &&\mathcal{L}^* = (\mathcal{G}^*_{\frac{1}{n}} \mid \Gamma^*,A \Rightarrow_{\frac{1}{n}} B,\Delta^*\mid \Gamma^* \Rightarrow_{\frac{1}{n}} \Delta^*)(\bar t_i)
    \end{eqnarray*} 
By well-formedness again, $\mathbb{X}$ contains some unique element $\mathfrak{L}^*\mid (\mathcal{R}^*)^{\tau^*}$ such that 
    \begin{eqnarray*}
     &&\tau^* =\eta^* \mathbin{\frown} (\mathcal{G}^*\mid \Gamma^*, A\to B\Rightarrow \Delta^*)\mathbin{\frown}(\mathcal{G}^* \mid \Gamma^* \Rightarrow \Delta^*)\\
    &&\mathcal{R}^* = (\mathcal{G}^*_{\frac{1}{n}} \mid \Gamma^* \Rightarrow_{\frac{1}{n}} \Delta^*)(\bar t_i)
    \end{eqnarray*} 
We similarly derive the following 
$$\mathfrak{L}^*\mid ((\mathcal{G}^*_{\frac{1}{n}} \mid \Gamma^* \Rightarrow_{\frac{1}{n}} A\to B,\Delta^*)(\bar t_i))^{\eta^* \mathbin{\frown} (\mathcal{G}^*\mid \Gamma^*\Rightarrow  A\to B, \Delta^*)}$$
and replace both $\mathfrak{L}^*\mid (\mathcal{L}^*)^{\sigma^*}$ and $\mathfrak{L}^*\mid (\mathcal{R}^*)^{\tau^*}$ in $\mathbb{X}$ with it.

    \item $\mathfrak{R}$ is of the form  
    $\mathfrak{L}^*\mid (\mathcal{L}^*)^{\sigma^*}$ where 
\begin{eqnarray*}
    &&\sigma^* =\eta^* \mathbin{\frown} (\mathcal{G}^*\mid \Gamma^*\Rightarrow \Delta^*)\mathbin{\frown}(\mathcal{G}^* \mid \Gamma^*\Rightarrow \Delta^*)\\
   &&\mathcal{L}^* = (\mathcal{G}^*_{\frac{1}{n}} \mid \Gamma^* \Rightarrow_{\frac{1}{n}} \Delta^*)(\bar t_i)
    \end{eqnarray*}
We replace $\mathcal{R}$ with the following
$$\mathfrak{L}^*\mid ((\mathcal{G}^*_{\frac{1}{n}} \mid \Gamma^*\Rightarrow_{\frac{1}{n}} \Delta^*)(\bar t_i))^{\eta^* \mathbin{\frown} (\mathcal{G}^*\mid \Gamma^*\Rightarrow \Delta^*)}$$
        \end{itemize}

Define $F(\mathbb{X})$ to be the result of substituting $\mathfrak{G}$ and other elements in $\mathbb{X}$. Now we check well-formedness of $F(\mathbb{X})$:
\begin{itemize}
    \item Property (1) is clear as the replacement of $\mathfrak{G}$ is still a left-most element. 
    \item The only way Property (2) fails is when we remove labels on the right branch at a split while keeping the left branch. Since we simultaneously remove both branching labels, Property (2) for $F(\mathbb{X})$ is preserved. 
    \item Elements in $F(\mathbb{X})$ are kept with the same number of components and the same stage. 
    \item Elements in $F(\mathbb{X})$ are all provided with a derivation.
\end{itemize}

% Left \exists
\textbf{(4):} If conditions in $(1)$, $(2)$, and $(3)$ fail, then components of $\mathfrak{G}$ must have labels either $(\mathcal{H})$ or of the form
    $$\eta \mathbin{\frown} (\mathcal{G}\mid \Gamma,\exists y A\Rightarrow \Delta)\mathbin{\frown}(\mathcal{G} \mid \Gamma,A[f(\bar x)/y]\Rightarrow \Delta)$$
Assume not all labels are $(\mathcal{H})$ and say $\mathfrak{G}$ is of the form
$$\mathcal{L}_1^{\sigma_1}\mid \cdots \mid\mathcal{L}_k^{\sigma_k}$$ and let $f_i(\bar x_i)$ be the corresponding Skolem function (add a dummy function in the case of $(\mathcal{H})$ label). Then each $\mathcal{L}_i$ is either $\mathcal{H}_{\frac{1}{n}}$ or of the form
$$((\mathcal{G}_{i})_\frac{1}{n} \mid \Gamma_i ,A_i[f_i(\bar x_i)/y_i]\Rightarrow_{\frac{1}{n}} \Delta_i)[\bar t_i/\bar x]
$$

By construction, $\bar x_i$ enumerates all free variables in $\mathcal{L}_i$ before substitution and thus $\bar t_i$ enumerates all unbounded terms in $\mathcal{L}_i$ except $f_i(\bar x_i)[\bar t_i/\bar x]$. Hence,  $f_i(\bar x_i)[\bar t_i/\bar x]$ does not appear elsewhere in $\mathcal{L}_i$ otherwise $f_i(\bar x_i)[\bar t_i/\bar x]$ would be a subterm of $\bar t_i$.

We now take $m$ such that $f_m(\bar x_m)[\bar t_m/\bar x]$ is not a proper subterm of any $f_i(\bar x_i)[\bar t_i/\bar x]$. If there exists $j\neq m$ such that 
$$f_m(\bar x_m)[\bar t_m/\bar x]=f_j(\bar x_j)[\bar t_j/\bar x]$$ then we must have $\sigma_m=\sigma_j$ and thus $\mathcal{L}_m=\mathcal{L}_j$ since the witness $\bar t_m$ and $\bar t_j$ are identified at the label $\sigma_m$. Let $j_1,\ldots,j_l$ enumerate all such $j\neq m$ and write 
$$\mathfrak{G}=\mathfrak{L}\mid \mathcal{L}_m^{\sigma_m}\mid  \mathcal{L}_{j_1}^{\sigma_{j_1}}\mid \cdots \mid \mathcal{L}_{j_l}^{\sigma_{j_l}}$$

We remove all such identical labeled hypersequents from $\mathfrak{G}$ by using $ec$ since $\mathcal{L}_m=\mathcal{L}_{j_1}=\cdots =\mathcal{L}_{j_l}$.
\begin{prooftree}
    \AxiomC{$\mathfrak{L}\mid \mathcal{L}_m\mid  \mathcal{L}_{j_1}\mid \cdots \mid \mathcal{L}_{j_l}$}
    \RightLabel{$ec$}
    \UnaryInfC{$\mathfrak{L}\mid \mathcal{L}_m$}
\end{prooftree}

Now we have that $f_m(\bar x_m)[\bar t_m/\bar x]$ does not appear elsewhere in the entire $\mathfrak{L}\mid \mathcal{L}_m$. We can introduce $\exists y_m$ by using Lemma \ref{ex1n} and add back $\mathcal{L}_{j_1},\ldots,\mathcal{L}_{j_l}$ by $ew$.
\begin{prooftree}
    \AxiomC{$\mathfrak{L}\mid ((\mathcal{G}_m)_{\frac{1}{n}} \mid \Gamma_m,A_m[f(\bar x_m)/y_m],\Rightarrow_{\frac{1}{n}} \Delta_m)(\bar t_m)$}
    \RightLabel{$ \exists\Rightarrow_{\frac{1}{n}}$}
    \UnaryInfC{$\mathfrak{L}\mid ((\mathcal{G}_m)_{\frac{1}{n}} \mid \Gamma_m,\exists y_m A_m,\Rightarrow_{\frac{1}{n}} \Delta_m)(\bar t_m)$}
    \RightLabel{$ew$}
    \UnaryInfC{$\mathfrak{L}\mid ((\mathcal{G}_m)_{\frac{1}{n}} \mid \Gamma_m,\exists y_m A_m,\Rightarrow_{\frac{1}{n}} \Delta_m)(\bar t_m)\mid  \mathcal{L}_{j_1}\mid \cdots \mid \mathcal{L}_{j_l}$}
\end{prooftree}
In comparison to the initial label $\sigma_m$ in $\mathfrak{G}$ $$\sigma_m=\eta_m \mathbin{\frown} (\mathcal{G}_m\mid \Gamma_m,\exists y_m A_m\Rightarrow \Delta_m)\mathbin{\frown} (\mathcal{G}_m\mid \Gamma_m, A_m[f_m(\bar x_m)/y_m]\Rightarrow \Delta_m)$$
we replace $\mathfrak{G}$ with 
$$\mathfrak{L}\mid (((\mathcal{G}_m)_{\frac{1}{n}} \mid \Gamma_m,\exists y_m A_m,\Rightarrow_{\frac{1}{n}} \Delta_m)(\bar t_m))^{\eta_m \mathbin{\frown} (\mathcal{G}_m\mid \Gamma_m,\exists y_m A_m\Rightarrow \Delta_m)}\mid  \mathcal{L}_{j_1}^{\sigma_{j_1}}\mid \cdots \mid \mathcal{L}_{j_l}^{\sigma_{j_l}}$$

For other elements $\mathfrak{R}$ in $\mathbb{X}$, there are two cases. 
    \begin{itemize}
    \item $\mathfrak{R}$ is of the form  
    $\mathfrak{L}^*\mid (\mathcal{L}_m^*)^{{\sigma_m}^*}$ where $\sigma_m^*$ has the form
    $$\eta^*_m \mathbin{\frown} (\mathcal{G}^*_m\mid \Gamma_m^*,\exists y_m A_m\Rightarrow \Delta_m^*)\mathbin{\frown} (\mathcal{G}_m^*\mid \Gamma_m^*, A_m[f_m(\bar x_m)/y_m]\Rightarrow \Delta_m^*)$$
Applying $ec$ to components $j$ that are identical to component $m$ and notice that then $f_m(\bar x_m)[\bar t_m/\bar x]$ now still satisfies the eigen-condition because $\mathfrak{R}$ does not contains more terms than $\mathfrak{G}$ due to Lemma \ref{sync}. Hence, we can apply $\exists\Rightarrow_{\frac{1}{n}}$ and $ew$ to obtain 
$$\mathfrak{L}^*\mid (((\mathcal{G}^*_m)_{\frac{1}{n}} \mid \Gamma^*_m,\exists y_m A_m \Rightarrow_{\frac{1}{n}} \Delta_m^*)(\bar t_m))^{\eta_m^* \mathbin{\frown} (\mathcal{G}_m^*\mid \Gamma_m^*, \exists y_m A_m\Rightarrow \Delta_m^*)}$$
and replace $\mathfrak{R}$ with it.

    \item $\mathfrak{R}$ is of the form  
    $\mathfrak{L}^*\mid (\mathcal{L}_m^*)^{{\sigma_m}^*}$ where $\sigma_m^*$ has the form
    $$\eta^*_m \mathbin{\frown} (\mathcal{G}^*_m\mid \Gamma_m^*\Rightarrow \Delta_m^*)\mathbin{\frown} (\mathcal{G}_m^*\mid \Gamma_m^*\Rightarrow \Delta_m^*)$$
We replace $\mathfrak{R}$ with the following
$$\mathfrak{L}^*\mid (((\mathcal{G}^*_m)_{\frac{1}{n}} \mid \Gamma^*_m \Rightarrow_{\frac{1}{n}} \Delta^*)(\bar t_i))^{\eta_m^* \mathbin{\frown} (\mathcal{G}_m^*\mid \Gamma_m^*\Rightarrow \Delta_m^*)}$$
        \end{itemize}

Define $F(\mathbb{X})$ to be the result of substituting $\mathfrak{G}$ and other elements in $\mathbb{X}$. Similar to the case in (1) and (2), $F(\mathbb{X})$ remains to be well-formed. 

% Trivial
\textbf{(5):} If all previous conditions fail, then $F(\mathbb{X})=\mathbb{X}$. Notice that by the construction of $T_{\mathcal{H}}$, the left-most element $\mathcal{G}$ in $\mathbb{X}$ must have all its labels to be $(\mathcal{H})$ in this case.

\end{definition}

\begin{lemma}\label{rec}
Given any well-formed set $\mathbb{X}$, either $\mu(F(\mathbb{X}))<\mu(\mathbb{X})$ or the left-most element in $\mathbb{X}$ only has labels taking form of $(\mathcal{H})$.
\begin{proof}
Proof by induction on the construction of $F$, since $F$ always reduce labels except in case (5) where the left-most element in $\mathbb{X}$ only has labels as $(\mathcal{H})$.  
\end{proof}
\end{lemma}

\begin{theorem}[Completeness of GŁ$\forall$]
If $\mathcal{H}$ is a valid first-order hypersequent sequent, then $ \mathcal{H}_{\frac{1}{n}}$ is derivable in GŁ$\forall$ for all $n$.

\begin{proof}
Let $\mathcal{H}$ and $n$ be given. We consider its tree $T_{\mathcal{H}}$ with free variables $\bar x$ and its leaf set $X_{h}$. By Lemma \ref{apsa} and Lemma \ref{aphb} applying to $T_{\mathcal{H}}$ and $X_{h}$, there exist terms $\bar t_1, \ldots,\bar t_k$ such that 

$$\mathcal{G}^1_{\frac{1}{n}}(\bar t_1)\mid\cdots  \mid \mathcal{G}^k_{\frac{1}{n}}(\bar t_k)$$
is valid hypersequent for any $G^1,\ldots, G^k\in X_{h}$.

Define $\mathbb{X}_0$ as in Definition \ref{init} and it is well-formed by Lemma \ref{inwf}. We now apply $F$ as in Definition \ref{lbrd} recursively to $\mathbb{X}_0$. By Lemma \ref{rec}, there exists some finite $t$ such that the left-most element $\mathfrak{G}$ of $F^t(\mathbb{X})$ has only labels as $(\mathcal{H})$. By well-formedness of $F^t(\mathbb{X})$, we have a derivation of $$\mathfrak{G}=(\mathcal{H}_{\frac{1}{n}}(\bar t_1))^{(\mathcal{H})}\mid\cdots  \mid (\mathcal{H}_{\frac{1}{n}}(\bar t_k))^{(\mathcal{H})}=\mathcal{H}_{\frac{1}{n}}^{(\mathcal{H})}\mid\cdots  \mid \mathcal{H}_{\frac{1}{n}}^{(\mathcal{H})}$$
since $\mathcal{H}$ has no free variables. We then derive $\mathcal{H}_{\frac{1}{n}}$ by $ec$.
\begin{prooftree}
    \AxiomC{$\mathcal{H}_{\frac{1}{n}}\mid\cdots  \mid \mathcal{H}_{\frac{1}{n}}$}
    \RightLabel{$ec$}
    \UnaryInfC{$\mathcal{H}_{\frac{1}{n}}$}
\end{prooftree}

\end{proof}
\end{theorem}

The main result Theorem 10 in \cite{main} can be shown by noticing $\Rightarrow A\oplus A^n$ is derivable from $\Rightarrow_{\frac{1}{n}}A$.
\begin{corollary}
The followings are equivalent:
\begin{enumerate}
    \item $A$ is a valid formula in $[0,1]$-semantics.
    \item $\Rightarrow_{\frac{1}{n}}A$ is derivable in GŁ$\forall$ for all $n$. 
    \item $\Rightarrow A\oplus A^n$ is derivable in GŁ$\forall$ for all $n$.
\end{enumerate}
\end{corollary}

\begin{corollary}
GŁ$\forall$ together with the following infinitary rule

\begin{prooftree}
    \AxiomC{$\mathcal{H}_{\frac{1}{n}}$ for all $n$}
    \RightLabel{$apx$}
    \UnaryInfC{$\mathcal{H}$}
\end{prooftree}

is complete with respect to $[0,1]$-semantics.
\end{corollary}
\begin{corollary}
Infinitary GŁ$\forall$ admits $cut$. 
\begin{proof}
The $cut$ rule is sound with respect to the semantics.
\end{proof}
\end{corollary}

\section{Further Directions}
We end this paper by mentioning a few potential research directions. 

First, continuous logic (see \cite{clcm}) is based on first-order Łukasiewicz logic, plus a $\frac{1}{2}$-operator, a built-in bounded metric $d$, and uniformly continuity requirement for all function and relation symbols. A generalization of GŁ$\forall$ to continuous logic and a similar proof of approximate completeness theorem should be possible. 

Second, it should be possible to strengthen the cut admissibility to a cut elimination result. We might need an effective proof of Theorem \ref{aphb} by the means of functional interpretation (see \cite{pmht}).

Lastly, GŁ and GŁ$\forall$ are non-constructive logics in the sense that $(A\to B)\vee (B\to A)$ is derivable due to the use of $split$. One might consider a ``intuitionistic" Łukasiewicz logic by removing $split$ and only allow derivations of single sequents. A development of the intutionistic proof system and semantics, together with some equivalent form of realizability and double-negation translations could be interesting and have potentials in proof mining for continuous model theory. Notice that a system of intuitionistic Łukasiewicz logic was already investigated by Arthan and Oliva in \cite{inlk} by removing double negation elimination. Their approach is different and it is also interesting to compare the two approaches. 

\section{Acknowledgment}
I would like to express my gratitude to Professor Henry Towsner for his invaluable supervision and support throughout this project. I am also deeply thankful to Eben Blaisdell for his insightful guidance and encouragement.

\end{document}